\documentclass[sn-mathphys-num]{sn-jnl}

\usepackage[T1]{fontenc}
\usepackage{amssymb,amsmath,latexsym}
\usepackage{amsthm}
\usepackage{mathrsfs}
\usepackage{diagbox}
\usepackage[shortlabels]{enumitem}
\usepackage{color}
\usepackage{comment}
\usepackage{anyfontsize}
\usepackage{xcolor}
\usepackage{stackengine}
\usepackage{calc}
\usepackage{physics}
\usepackage{accents}
\usepackage{hyperref}
\usepackage{appendix}
\usepackage{colortbl}
\usepackage{booktabs}
\usepackage{tabularx}
\usepackage{tikz-cd}
\usepackage{aliascnt}
\usepackage{cleveref}
\usetikzlibrary{arrows.meta,positioning}
\crefname{appendix}{appendix}{appendices}
\Crefname{appendix}{Appendix}{Appendices}

\input{shorthand.sty}
\begin{document}
\title[Generalized Poincar\'e inequality for quantum Markov semigroups]{Generalized Poincar\'e inequality for quantum Markov semigroups}
\author[1]{\fnm{Marius} \sur{Junge}}\email{mjunge@illinois.edu}
\author*[2]{\fnm{Jia} \sur{Wang}}\email{jiawang5@illinois.edu}

\affil[1]{\orgdiv{Department of Mathematics}, \orgname{University of Illinois, Urbana-Champaign}}
\affil[2]{\orgdiv{Department of Physics}, \orgname{University of Illinois, Urbana-Champaign}}
\footnotetext[1]{The authors are listed in alphabetical order, contributed equally to this work, and should be considered co-first authors.}

\abstract{
We prove a noncommutative $p$-Poincar\'e inequality for GNS-detailed-balance quantum
Markov semigroups (QMSs) on non-tracial $\sigma$-finite von Neumann algebras, assuming only the existence of a spectral gap.
Extending the semi-commutative results of Huang and Tropp \cite{huang2021poincare}, we first establish the inequality for tracial von Neumann algebras.
We use Markov dilations to obtain chain-rule estimates for Dirichlet forms and amalgamated free products to define an appropriate noncommutative derivation.
We then extend the result to QMSs satisfying GNS detailed balance on non-tracial $\sigma$-finite von Neumann algebras, using Haagerup's reduction and Kosaki's interpolation theorem.
As applications, we recover sub-exponential concentration inequalities and estimate the Lipschitz and completely bounded Lipschitz diameters of the QMSs.
}


\maketitle
\setcounter{tocdepth}{2}
{\small\tableofcontents}
\section{Introduction}\label{sec:intro}
$L^p$ Poincar\'e inequalities have a long history in harmonic analysis, geometry, and probability.
They control higher moments and yield concentration results.  Building on
Lust-Piquard's dimension-free estimates for discrete Riesz transforms
\cite{lustpiquard1998riesz,lustpiquard2004dimensionfree}, obtained
using noncommutative techniques by Pisier
\cite{pisier1988riesz}, Ben Efraim and
Lust-Piquard \cite{benefraim-lustpiquard2007} proved dimension-free $L^p$
Poincar\'e inequalities for the discrete cube and corresponding inequalities for
the CAR algebra.  These arguments adapt Gaussian methods of Maurey and Pisier to special structures, such as
Bernoulli and fermionic settings. Similarly, Junge and Zeng \cite{junge2015noncommutative} established broader noncommutative results that rely on a
positive curvature condition, which is stronger than a spectral gap. Huang and Tropp \cite{huang2021poincare} later obtained
semi-commutative matrix-valued $L^p$ moment bounds from a spectral gap alone.
In this work, we extend this spectral-gap principle from the semi-commutative setting to fully
noncommutative quantum Markov semigroups satisfying GNS detailed balance on
non-tracial $\sigma$-finite von Neumann algebras.

We start with the Poincar\'e spectral-gap estimate \cite{poincare1890}.
In its classical Euclidean form, the Poincar\'e inequality states that if
$\Omega\subseteq\mb R^d$ is a bounded, connected open set with Lipschitz boundary,
then there exists a constant $\alpha>0$, depending only on $\Omega$, such that every
smooth compactly supported function $f$ on $\Omega$ satisfies
\[\alpha\int_{\Omega}\bigl|f(x)-f_{\Omega}\bigr|^2\,dx
  \le \int_{\Omega}\|\nabla f(x)\|_2^2\,dx,\]
where $f_{\Omega}=\frac{1}{|\Omega|}\int_{\Omega}f(x)\,dx$ denotes the mean of $f$.

Meyer \cite{meyer1976carre} reformulated the classical gradient in the
semigroup setting as the \textit{carr\'e du champ}, which we refer to as the gradient
form. This framework was subsequently developed systematically by Bakry and
{\'E}mery \cite{bakry-emery}. Let
$(T_t)_{t\ge0}=e^{-tL}$ be a conservative, $\mu$-symmetric Markov semigroup
on a probability space $(X,\mu)$, and let $\Gamma$ be its gradient form.
The semigroup is said to have spectral gap $\alpha>0$ if
\begin{equation}\label{eq:measure-Poincare}
  \|f-\mu(f)\|_{L^2(X,\mu)}
  \le \frac{1}{\sqrt{\alpha}}
  \|\Gamma(f,f)\|_{L^1(X,\mu)}^{1/2}.
\end{equation}
More generally, for $p,q/2\in[1,\infty]$, we say that $T_t$ satisfies the
\emph{Poincar\'e inequality} $\Poin(p,q)$ with constant $C>0$ if
\begin{equation}\label{eq:commutative-P-pq}
  \|f-\mu(f)\|_{L^p(X,\mu)}
  \le C\|\Gamma(f,f)\|_{L^{q/2}(X,\mu)}^{1/2}.
\end{equation}
We use the same terminology in the noncommutative setting.  Let
$(\M,\tau)$ be a tracial von Neumann algebra with faithful normal tracial
state $\tau$, and let $(T_t)_{t\ge0}=e^{-tL}$ be a weak$^*$-continuous,
$\tau$-symmetric quantum Markov semigroup.  Its gradient form is
\[
  \Gamma(x,y)
  =\tfrac12\bigl(L(x^*)y+x^*L(y)-L(x^*y)\bigr),
  \qquad x,y\in\Dom(L).
\]
The fixed-point space
\[
  \N=\{x\in\M:T_t(x)=x\text{ for all }t\ge0\}
\]
is a von Neumann subalgebra. We denote by $E:\M\to\N$ the canonical
$\tau$-preserving conditional expectation.

We have the following noncommutative analogue of \cref{eq:commutative-P-pq}.
\begin{definition}[Tracial Poincar\'e inequality $\Poin(p,q)$]\label{def:Poin(p,q)}
For $p,q/2\in[1,\infty]$, the QMS satisfies the \emph{Poincar\'e
inequality} $\Poin(p,q)$ with constant $C>0$ if, for every self-adjoint
$x\in\Dom(L)$,
\begin{equation}\label{eq:P-pq}
  \|x-E(x)\|_{L^p(\M,\tau)}
  \le C\|\Gamma(x,x)\|_{L^{q/2}(\M,\tau)}^{1/2}.
\end{equation}
In particular, the spectral gap is exactly $\Poin(2,2)$ (with
$C=\alpha^{-1/2}$).
\end{definition}
With the above definition, the tracial form of our main result is the following:
\medskip
\begin{theorem}[$\Poin(p,p)$]\label{main:tracial}
Let $(\M,\tau)$ be a tracial von Neumann algebra and $(T_t)_{t\ge0}=e^{-tL}$ a $\tau$--symmetric quantum Markov semigroup with fixed-point algebra $\N$ and conditional expectation $E:\M\to\N$.
If $\{T_t\}_{t\ge0}$ has a spectral gap $\alpha>0$, then for every self-adjoint $x$ in the domain of $L$ and $p=2$ or $p\ge3$, it satisfies the Poincar\'e inequality $\Poin(p,p)$ with constant $\tfrac{p}{\sqrt{2\alpha}}$, i.e.,
\[
\|x-E(x)\|_{L^p(\M,\tau)}  \le  \frac{p}{\sqrt{2\alpha}} \|\Gamma(x,x)\|^{1/2}_{L^{p/2}(\M,\tau)}.
\]
\end{theorem}

The proof uses the argument of Huang and Tropp as a stepping stone. In particular, they established a chain rule estimate relating the Dirichlet form to the associated gradient form. A further ingredient in their argument is a symmetrization procedure, based on doubling variables in the product space of a Markov process, which leads to key cancellations.
Our proof, in addition, identifies the appropriate dilation framework for quantum Markov semigroups, which allows us to derive chain rule estimates directly in the noncommutative setting. Since there are obstacles to generalizing their product measure symmetrization procedure, we instead formulate it in terms of amalgamated free products, providing a natural setting for noncommutative algebras.

\paragraph{Beyond the tracial setting: GNS-detailed balance via Haagerup reduction.}
Many physically relevant QMS are not tracial but are reversible with respect to a faithful normal state $\phi$ (GNS-symmetry/detailed balance). In that case, one works on Haagerup $L^p$ spaces~\cite{haagerup1979lpspaces} $L^p(\M)$ and uses the GNS inner product in Section~\ref{sec:lp-spaces}.

Our strategy reduces the non-tracial problem to the tracial one via Haagerup's
reduction (crossed-product) technique.
We embed $(\M,\phi)$ into a semifinite crossed product $(\widehat{\M},\widehat{\phi})$,
which admits an increasing sequence of tracial subalgebras $(\M_k,\tau_k)$.
The semigroup is lifted to be $\widehat{\phi}$-detailed balance on $\widehat{\M}$ and $\tau_k$-symmetric on $\M_k$,
where a $\Poin(p,q)$ inequality on $\wh\M$ is obtained by tracial $\Poin(p,p)$
on each slice $(\M_k,\tau_k)$.
The result is then descended to $\M$ via conditional expectations compatible
with the modular action, giving a GNS-$\phi$-detailed balance version of $\Poin(p,p)$.
Throughout the paper, we abbreviate the associated Haagerup gradient seminorm by
\[
  \|a\|_{\Gamma_p}
  :=\|\Gamma_p(a,a)\|_{L^{p/2}(\M)}^{1/2}.
\]
Thus the inequality takes the form
\[
  \|a-E_p^\phi(a)\|_{L^p(\M)}
  \le C\|a\|_{\Gamma_p},
\]
where $E_p^\phi$ denotes the $L^p$ extension of the $\phi$-preserving conditional expectation onto the fixed-point algebra,
$\Gamma_p$ the $L^p$ version of the gradient form, and $\|\cdot\|_{L^p(\M)}$ the norm in Haagerup's $L^p$ space.
A key advantage of Haagerup's $L^p$ construction is its functoriality: every normal $*$-homomorphism between von Neumann algebras induces canonically bounded maps on the associated $L^p$-spaces, compatible with composition, duality, and interpolation.
Our main result is stated in this framework as follows:
\medskip
\begin{theorem}[cf \Cref{thm:GNS-Lp-ppp} and \Cref{cor:non-self-adjoint}]\label{main:GNS}
Let $(\M,\phi)$ be a $\sigma$-finite von Neumann algebra with faithful normal state $\phi$, and let $(T_t)_{t\ge0}$ be a GNS-$\phi$-detailed-balance QMS with fixed-point algebra $\N$. Suppose $(T_t)$ has an $L^2(\M,\phi)$ spectral gap $\alpha>0$. Let $E_p^\phi:L^p(\M)\to L^p(\N)$ be the induced conditional expectation. Then for every $p=2$ and every $p\ge3$, and for all self-adjoint $a\in L^p(\M)_{\mathrm{sa}}$ in the domain, we have
\[
\|a-E_p^\phi(a)\|_{L^p(\M)} \leq \frac{p}{\sqrt{2\alpha}}\|a\|_{\Gamma_p}.
\]
Furthermore, for general $a\in L^p(\M)$ in the domain,
\[
\|a-E_p^\phi(a)\|_{L^p(\M)} \leq \frac{p}{\sqrt{2\alpha}}
\bigl(\|a\|_{\Gamma_p}+\|a^*\|_{\Gamma_p}\bigr).
\]
\end{theorem}

\paragraph{Sub-exponential concentration inequalities.}
Poincar\'e-type inequalities are closely connected to the concentration of
measure. In the semi-commutative setting, such a principle was established by
Huang and Tropp \cite[Thm~2.7]{huang2021poincare} for Markov processes with a spectral
gap, showing sub-exponential concentration.
Their approach relies on matrix inequalities and dilation techniques, and improves
upon earlier results of Aoun, Banna, and Youssef \cite{Aoun_Banna_Youssef_2020}.
In particular, Huang and Tropp introduce a symmetrization (or doubling-of-variables)
argument that avoids technical difficulties present in the earlier approach.

We show that the $\Poin(p,p)$ inequality proved in this work implies a
\emph{sub-exponential} concentration inequality for observables measured in the
Lipschitz seminorm induced by the gradient form as below.
\begin{cor*}[cf \Cref{cor:concentration}]
Let $T_t=e^{-tL}$ be a GNS--$\phi$-detailed balance QMS with spectral gap $\alpha>0$ and $\phi$-preserving conditional expectation $E^\phi$ onto its fixed-point algebra. Then for any $x\in\M$ and large $t>0$,
\[
\bP_{\phi}(|x-E^\phi(x)|>t)
\le
2\exp\Bigl(-\frac{\sqrt{\alpha} t}{8\sqrt{2}e \|x\|_{\Lip_\Gamma}}\Bigr).
\]
\end{cor*}
Here the Lipschitz seminorm is defined by
\[
\|x\|_{\Lip_\Gamma} = \max\Bigl\{\|\Gamma(x,x)\|_\infty^{1/2}, \|\Gamma(x^*,x^*)\|_\infty^{1/2}\Bigr\},
\]
to be distinguished from the Lipschitz complexity defined in \cite[eq.~(1.1)]{Araiza2026resourcedependent}.

This reveals a qualitative distinction between spectral-gap and stronger
functional-inequality regimes.
An $L^2$ spectral gap alone does not yield Gaussian or Talagrand-type concentration in general, while one can still obtain sub-exponential tails.
We compare this regime with the stronger concentration implied by modified logarithmic Sobolev inequalities (MLSI) and discuss a birth--death process \Cref{ex:bd} that rules out Gaussian decay under a uniform spectral gap assumption alone.

For finite-dimensional matrix algebras, we establish estimates for the Lipschitz diameters and completely bounded Lipschitz diameters \Cref{cor:lip-cb-diam}.
We summarize the main results and local implications in \Cref{fig:road}.
\begin{figure*}[!h]
\centering
\begin{tikzpicture}[
    scale=.88,
  transform shape,
  box/.style={
    draw=black!75,
    rounded corners=3pt,
    line width=.7pt,
    align=center,
    inner xsep=5pt,
    inner ysep=5pt,
    minimum height=13mm,
    font=\normalfont\scriptsize
  },
  assumption/.style={
    box,
    fill=black!3,
    text width=22mm
  },
  result/.style={
    box,
    fill=black!8,
    text width=24mm
  },
  application/.style={
    box,
    fill=black!5,
    text width=25mm
  },
  edge label/.style={
    fill=white,
    inner sep=2pt,
    align=center,
    font=\normalfont\footnotesize
  },
  note/.style={
    draw=black!45,
    densely dotted,
    rounded corners=2pt,
    align=center,
    inner sep=4pt,
    font=\normalfont\scriptsize,
    text width=42mm
  },
  implies/.style={
    -{Latex[length=2.2mm,width=1.6mm]},
    line width=.8pt
  },
  nonimp/.style={
    densely dashed,
    line width=.8pt
  }
]

\node[assumption] (gap)
  {Spectral gap\\[-1pt]$\alpha>0$};

\node[result,right=15mm of gap] (pp)
  {$O(p)\,\Poin(p,p)$\\[-1pt]
   \Cref{main:GNS}};

\node[result,right=4mm of pp] (pinf)
  {$O(p)\,\Poin(p,\infty)$\\[-1pt]
   \Cref{cor:lip-norm-estimate}};

\node[application,right=4mm of pinf] (conc)
  {Sub-exponential\\
   concentration\\[-1pt]
   \Cref{cor:concentration}};

\draw[implies] (gap) --
  node[edge label,above=1mm]
  {GNS QMS}
  (pp);

\draw[implies] (pp) -- (pinf);
\draw[implies] (pinf) -- (conc);

\node[application,below=10mm of conc] (diam)
  {Lipschitz diameter and\\
   completely bounded\\
   Lipschitz diameter\\[-1pt]
   \Cref{cor:lip-cb-diam}};

\draw[implies]
  (pinf.south east) --
  node[edge label,midway,xshift=10mm]
    {finite\\dimension}
  (diam.north west);

\node[application,below right=13mm and 18mm of gap] (tal)
  {Talagrand\\
   inequality\\[-1pt]
   \Cref{def:talagrand}};

\draw[nonimp]
  (pinf) --
  node[edge label,midway,xshift=-10mm]
    {Not equivalent\\
     by \Cref{ex:bd}}
  (tal.north east);

\node[assumption,below=30mm of gap] (mlsi)
  {$\alpha_{\mathrm{MLSI}}>0$\\\cite{Gao_Junge_LaRacuente_Li_2025}};

\node[result,right=8mm of mlsi] (sqrtp)
  {$O(\sqrt p)\,\Poin(p,\infty)$\\[-1pt]
   \cite[Thm~5.12]{Gao_Junge_LaRacuente_Li_2025}};

\draw[implies] (mlsi) -- (sqrtp);
\draw[implies] (sqrtp.north) -- (tal.south);

\draw[implies]
  (sqrtp.east)
  to[out=0,in=-90,looseness=1.15]
  (pinf.south);
\end{tikzpicture}

\caption{Roadmap of the logical structure of the main results. Solid arrows denote
implications, while dashed lines indicate that no implication is
asserted. The birth--death example shows a violation of the dimension-free geometric Talagrand inequality, but the $O(p)~\Poin(p,p)$ estimate holds.}
\label{fig:road}
\end{figure*}
\paragraph{Related work.}
To our knowledge, this is the first result establishing a noncommutative $(p,p)$-type Poincar\'e inequality under only an $L^2$ spectral gap assumption, within the class of QMSs satisfying detailed balance. Earlier works obtain Poincar\'e inequalities with different index pairs and under stronger functional or curvature assumptions. In the noncommutative setting, such curvature bounds are often difficult to verify beyond specific model classes \cite{Wirth_Zhang_2023,münch2024intertwiningcurvatureboundsgraphs}.

In \cite{Gao_Junge_LaRacuente_Li_2025}, assuming a modified logarithmic Sobolev inequality (MLSI) with constant $\alpha_{MLSI}>0$, the authors prove $\PI(p,\infty)$ with optimal order $O(\sqrt{p})$.
The relation between our $\PI(p,p)$ and the MLSI-based $\PI(p,\infty)$ regime is not clear, and neither inequality evidently forces the other. MLSI further implies the geometric Talagrand inequality from \cite[Remark 6.9]{gao_fisher_2020}, while we illustrate in \Cref{ex:bd} that Talagrand-type bounds do not hold for all QMS with only an $L^2$ spectral gap.

Under the noncommutative $\Gamma_2$-condition
\[
  \tau\bigl(\Gamma(T_t x,T_t x) y\bigr) \le  C e^{-2\alpha t}
  \tau\bigl(T_t \Gamma(x,x) y\bigr)\qquad (x\in\Dom(A^{1/2}),  y\in\N_+),
\]
\cite{junge2015noncommutative} obtains both $(p,p)$- and $(p,\infty)$-type inequalities with constants of order $O(p)$ and $O(\sqrt{p})$, respectively.
Since both MLSI and $\Gamma_2$ imply an $L^2$ spectral gap, while the converse fails in general, our hypothesis is strictly weaker.

As for applications, \cite{De_Palma_2023} uses a $(2,\infty)$-type Poincar\'e inequality to derive lower bounds on quantum circuit depth. Using our framework, we say that a $(2,\infty)$-type Poincar\'e inequality holds for a GNS-$\phi$-detailed-balance QMS if
\begin{align*}
    \|x - E^\phi(x)\|_{L^{2}(\mathcal M,\phi)_L} \le C \|\Gamma(x,x)\|^{1/2}_\infty,\quad   x = x^*,
\end{align*}
where $\|x\|_{L^{2}(\mathcal M,\phi)_L} = \sqrt{\tau(D_\phi x^*x)}$. We show $(2,\infty)$-type Poincar\'e inequality in \Cref{cor:lip-norm-estimate} by choosing $p = 2, a = x D_\phi^{1/2}$, where $D_\phi$ is the density of $\phi$. The relationship between the $(2,\infty)$-type Poincar\'e inequality used in \cite{De_Palma_2023} and our formulation \cref{eq:P-pq} (up to a dimension factor) is made explicit in \cite[(3.9), (4.16)]{Araiza2026resourcedependent}. Further connections to circuit complexity appear in \cite{Bu_2024,Li_2025,Ding_2025,araiza2025}.
When detailed balance (or symmetry) fails, $\mathrm{PI}(p,q)$ generally breaks down; instead, space--time Poincar\'e inequalities under weaker symmetry assumptions are established in \cite{li2025}.
We also mention related unpublished work of Li Gao, Zhendong Xu, and Hao Zhang,
communicated to the first-named author, which provides an independent proof for the tracial case.
\paragraph{Organization of this paper.}
\Cref{sec:preliminaries} collects notation and background on von Neumann algebras, quantum Markov semigroups, Markov dilations, and amalgamated free products.
\Cref{sec:ppp-tracial} introduces the swap operator and the amalgamated free product derivation, proves the generalized Klein inequality used in the argument, and establishes \Cref{main:tracial}.
\Cref{sec:GNS} reviews Haagerup and Kosaki $L^p$ spaces, states the GNS detailed balance framework, and proves \Cref{main:GNS} via Haagerup reduction.
\Cref{sec:examples} presents examples and applications, including
sub-exponential concentration and finite-dimensional estimates for the
Lipschitz diameter and its completely bounded analogue.
The appendices give the Markov-dilation chain estimate,
spectral-gap stability under regularization, domain convergence,
and the convexification argument needed in the general case.

\textbf{Notation.}
We use script letters $\M,\N$ for von Neumann algebras, $\tau$ for a trace and $\phi$ for a non-tracial state.
The identity operator is denoted by $1$ and the identity map by $\Id$.
For $x\in\M$, $x^*$ denotes the adjoint and $\Phi_*$ the pre-adjoint of a normal map $\Phi$.
Let $(T_t)_{t\ge0}=e^{-tL}$ be a QMS with generator
$L$,
and denote by $\Gamma$ and $\mathscr E$ the associated gradient and Dirichlet forms.

\section{Preliminaries}\label{sec:preliminaries}

\subsection{Von Neumann algebras and conditional expectations}\label{sec:von-neumann-algebras}
Let $(\mathcal{M},\phi)$ denote a von Neumann algebra $\mathcal{M}$ equipped with a normal, faithful state $\phi$.
Write $(\pi_\phi,H_\phi,\xi_\phi)$ for the GNS (Gelfand--Naimark--Segal) representation of $\M$ with respect to $\phi$, where $\pi_\phi:\M\to B(H_\phi)$ is a faithful normal $*$-representation and $\xi_\phi\in H_{\phi}$ is cyclic and separating for $\pi_\phi(\mathcal{M})$, satisfying
\[
  \phi(x)=\langle \xi_\phi,\pi_\phi(x) \xi_\phi\rangle ,\qquad x\in\mathcal{M}.
\]
Here $\xi_\phi$ is \emph{cyclic} for $\M$ if $\overline{\M\xi_\phi}=H_\phi$, and \emph{separating} if $x\xi_\phi=0$ implies $x=0$, or equivalently $\overline{\M'\xi_\phi}=H_\phi$.
After identifying $\mathcal{M}$ with $\pi_\phi(\mathcal{M})$, we may regard $\M$ as a von Neumann algebra on $H_\phi$ with cyclic and separating vector $\xi_\phi$.

Define the densely defined antilinear map
\[
  S: H_\phi\to H_\phi,\qquad S(x \xi_\phi)=x^* \xi_\phi,\qquad x\in\mathcal{M}.
\]
By Tomita--Takesaki theory \cite{takesaki_theory_2003-2}, the closure $\overline{S}$ admits the polar decomposition $\overline{S}=J \Delta^{1/2}$, where $J$ is a conjugate-linear unitary and $\Delta$ is a positive self-adjoint operator. We call
$J$ the \emph{modular conjugation} and $\Delta:=\overline{S}^*\overline{S}$ the \emph{modular operator} associated with $\phi$.
The modular automorphism group $(\sigma_t^\phi)_{t\in\R}$ is then given by
\begin{equation}\label{eq:modular-group}
      \sigma_t^\phi(x)=\Delta^{it} x \Delta^{-it},\qquad x\in\mathcal{M},  t\in\R.
\end{equation}

Let $\N \subset \M$ be a von Neumann subalgebra. If $\N$ is globally invariant under the modular automorphism group,
\[
  \sigma_t^\phi(\N)=\N \quad (t\in\R),
\]
then by Takesaki's theorem \cite{Takesaki_1972,takesaki_theory_2003-2} there exists a \emph{unique} {$\phi$--preserving conditional expectation}
\[
  E_\phi : \M \to \N
\]
such that
\begin{enumerate}[(i)]
  \item $E_\phi(y) = y$ for all $y \in \N$;
  \item $E_\phi(y_1 x y_2) = y_1 E_\phi(x) y_2$ for all $x \in \M$, $y_1,y_2 \in \N$;
  \item $\phi(E_\phi(x)) = \phi(x)$ for all $x \in \M$.
\end{enumerate}
It follows that $E_\phi$ is normal, completely positive, contractive, and $\N$--bimodular.
For more background on von Neumann algebras, see \cite{Takesaki_1979,takesaki_theory_2003-2}. In the remainder of this section, we discuss $L^p$ spaces in the tracial case, postponing the non-tracial case to \Cref{sec:lp-spaces}.

Let $(\M,\tau)$ be a von Neumann algebra equipped with a faithful normal tracial state $\tau$.
In this setting, one can define noncommutative $L^p$ spaces in complete analogy with the classical case.
For $1\le p<\infty$, the $L^p$ norm is
\[
  \|x\|_p := \tau(|x|^p)^{1/p},
  \qquad |x|=(x^*x)^{1/2}, \quad x\in\M,
\]
and $\|x\|_\infty := \|x\|$ is the operator norm.
The completion of $\M$ under this norm gives the noncommutative $L^p$ space $L^p(\M,\tau)$.

We write $\M_{\mathrm{sa}}$ for the self-adjoint part of $\M$ and $\M_+$ for the positive cone.
The case $p=2$ is of particular importance: $L^2(\M,\tau)$ carries the inner product
\[
  \langle x,y\rangle_{L^2(\M)} := \tau(x^*y), \qquad x,y\in\M,
\]
under which $\M$ embeds densely and isometrically.
This is exactly the GNS Hilbert space $H_\phi$ associated with $\tau$.

\subsection{Functional inequalities}
In this section, we introduce the Poincar\'e inequality for trace-symmetric quantum Markov semigroups.
Let $(T_t)_{t\ge0}:\M\to\M$ be a quantum Markov semigroup (QMS) on $\mathcal{M}$, i.e., a family of normal, unital, completely positive maps (UCP maps)
satisfying
\begin{enumerate}[i)]
  \item $T_0 = \mathrm{Id}_{\mathcal{M}}$ and $T_{t+s} = T_t\circ T_s$ for all $s,t\ge0$,
  \item for each $x\in\mathcal{M}$, the map $t\mapsto T_t(x)$ is weak$^*$--continuous.
\end{enumerate}
We write formally
$T_t  =  e^{-tL},$
where the semigroup generator (Lindbladian) $L$
is given by
\begin{equation}\label{eq:lindbladian}
    L(x)
 :=
\mathrm{w}^*\text{-}   \lim_{t\to0}\frac{1}{t}\bigl(x - T_t(x)\bigr),
\end{equation}
on the domain of those $x\in\mathcal{M}$ for which the limit exists.
Assume now that $\M$ is finite and is equipped with a faithful normal tracial state $\tau$; we write this pair as $(\M,\tau)$. We define a QMS to be $\tau$-symmetric if it is self-adjoint with respect to the trace, i.e.,
\begin{equation}\label{eq:tau-sym}
\begin{aligned}
    \tau(x^* T_t(y)) = \tau(T_t(x)^* y), \quad   x, y \in \mathcal{M}, t \ge 0.
\end{aligned}
\end{equation}

This definition implies that $L$ is positive and $\tau$-symmetric on $L^2(\mathcal{M},\tau)$. Its fixed-point algebra
\[
\mathcal{N} := \ker(L) = \{ x \in \mathcal{M} : T_t(x) = x \text{ for all } t \ge 0 \}
\]

is a von Neumann subalgebra of $\mathcal{M}$, and each $T_t$ is an $\mathcal{N}$--bimodule map,
\[
T_t(a x b) = a T_t(x) b, \quad   a, b \in \mathcal{N},\ x \in \mathcal{M}.
\]

In particular,
\[
T_t \circ E = E \circ T_t = E,
\]
where
$
E : \mathcal{M} \to \mathcal{N}.
$
is the $\tau$-preserving conditional expectation onto $\mathcal{N}$.
Before introducing the Poincar\'e inequality, we recall the associated energy and gradient forms.
\begin{definition}[Dirichlet form]\label{def:DirichletForm}
The Dirichlet form (or energy functional) of the semigroup is the quadratic map
\[
\mathscr E:\Dom(L^{1/2})\to\mathbb R_{\ge0},
\qquad
\mathscr E(x):=\|L^{1/2}x\|_{L^2(\mathcal M,\tau)}^2.
\]
For $x\in\Dom(L)$, this equals $\tau(x^*Lx)$.
\end{definition}
\begin{definition}[Gradient form]\label{def:Gamma}
The noncommutative gradient form associated to the Lindbladian $L$ is the sesquilinear map $\Gamma: \cA \times \cA \to \mathcal{M}$ where $
\mathcal{A}
:=
\Dom(L)\cap\mathcal{M},
$\footnote{$\A$ is a dense $*$--subalgebra of $\mathcal{M}$ and a core for $L^{1/2}$ \cite{cipriani_derivations_2003, davies_non-commutative_1992}.}
is defined by \cite[Section~2.1]{gao_fisher_2020} \[
\Gamma(x,y)
 =
\tfrac12\bigl(L(x^*)  y  +  x^*  L(y) - L(x^*y)\bigr),
\quad
  x,y\in\mathcal{A}.
\]
Equivalently, by the weak$^*$--limit
\begin{equation}\label{eq:Gamma}
    \Gamma(x,y)
 :=
\mathrm{w}^*\text{-}   \lim_{t\downarrow0}
\frac{1}{2t}\Bigl(
T_t(x^*y)
 -
T_t(x^*)  y
 -
x^*  T_t(y)+x^*y
\Bigr),
\quad
 x,y\in\mathcal{A},
\end{equation}
or equivalently from \cite[Lemma~9.1]{cipriani_derivations_2003},

\begin{equation}\label{eq:Gamma-epsilon}
    \Gamma(x,y)
 :=
\mathrm{w}^*\text{-}   \lim_{\epsilon\downarrow0}
\frac{1}{2}\Bigl(
x^*\frac{L}{1+\epsilon L}(y)+\frac{L}{1+\epsilon L}(x)^*y-\frac{L}{1+\epsilon L}(x^*y)
\Bigr),
\quad
  x,y\in\mathcal{A}.
\end{equation}

In particular, $\Gamma(x,x)\in\mathcal{M}_+$ for every $x\in\mathcal{A}$. We can also rewrite the energy form as
\[
\mathscr E(x) = \tau(\Gamma(x, x)),\quad x\in\A.
\]

\end{definition}

\begin{definition}
Let
$\{T_t\}_{t\ge0}$ a $\tau$--symmetric QMS with Lindbladian $ L$, gradient form $\Gamma$, and $E : \M \to \N$
be the $\tau$\nobreakdash-preserving conditional expectation onto $\N$.
Let $p,r/2\in[1,\infty]$,
we say that
\begin{itemize}
    \item $T_t$ admits the Poincar\'e inequality $\Poin (p,r)$ with constant $c>0$, if
    \[
        \| x - E(x) \|_{L^p(\M,\tau)} \le c   \| \Gamma(x, x) \|^{1/2}_{L^{r/2}(\M,\tau)},\quad  x\in\A
    \]
    \item $L$ admits the spectral gap $\alpha>0$ if
\begin{equation}
\alpha\|x - E(x)\|_{L^2(\M,\tau)}^{2}\leq \|\Gamma(x,x)\|_{L^1(\M,\tau)}.
\quad
 x\in\mathcal{A}.\label{eq:gap}
\end{equation}
\end{itemize}
\end{definition}
Note that the spectral gap is a special case of Poincar\'e inequalities $\PI(2,2)$, and the first main result we later derive is that given a spectral gap, one can derive a $\Poin(p,p)$ for $p=2$ and $p\geq 3$.
Integrating \Cref{eq:gap} gives the $L^2$ decay estimate
\[
\|T_t(x)-E(x)\|_{L^2(\M,\tau)} \le e^{-\alpha t}\|x-E(x)\|_{L^2(\M,\tau)},\qquad x\in\A,\ t\ge0.
\]

\subsection{Construction of amalgamated free product (AFP)}\label{sec:free-product}
\subsubsection{Motivation: Fock spaces}\label{sec:phys-motivation}
Starting from a one--particle Hilbert space $\H$ with vacuum vector $\xi$, one defines three standard Fock space constructions.
The symmetric Fock space arises as the GNS representation of the CCR (canonical commutation relations) algebra, the antisymmetric Fock space corresponds to the CAR (canonical anti--commutation relations) algebra, and the full (free) Fock space serves as the analogue without commutation restrictions.

\begin{itemize}[leftmargin=*]
  \item {Bosons (CCR).}
  \[
    \mathcal F_{\mathrm{Bose}}(\H) =  \bC \xi \oplus\bigoplus_{n=1}^\infty \mathrm{Sym}^n(\H).
  \]

  \item {Fermions (CAR).}
  \[
    \mathcal F_{\mathrm{Fermi}}(\H) =  \bC \xi \oplus\bigoplus_{n=1}^\infty \bigwedge^n \H.
  \]

\item {Free particles (free product algebra).}
\begin{equation}
    \label{eq:free-fock}
      \cF_{\mathrm{Free}}(\H)
  = \bC \xi \oplus \bigoplus_{n=1}^\infty \H^{\otimes n}.
\end{equation}
\end{itemize}

In the presence of multiple vacua encoded by a subalgebra, one passes to the
amalgamated free product (AFP) Fock space, which may be viewed heuristically
as obtained by gluing free Fock spaces over the different vacuum sectors.
\subsubsection{Amalgamated free products (AFP)}
\label{sec:afp}

We recall the construction of amalgamated free products (AFP) of tracial von Neumann algebras, following
\cite{speicher_combinatorial_1998,araki_operator_1985,ioana_amalgamated_2008,voiculescu_free_1992,brown_-algebras_2008}.
Let $(\M_{1},\tau_{1})$ and $(\M_{2},\tau_{2})$ be tracial von Neumann algebras with a common subalgebra
$\N\subset \M_{i}$, $i=1,2$, such that $\tau_{1}|_{\N}=\tau_{2}|_{\N}$.
There exist unique trace-preserving conditional expectations
$E_i\colon \M_i\to \N$ \cite[Theorem~4.2]{takesaki_theory_2003-2}, and we write
\[
  \mathring{\M}_i := \ker(E_i).
\]
Each $\M_i$ is naturally an $\N$--$\N$ bimodule with $\N$-valued inner product
$\langle x,y\rangle_{\N}=E_i(x^*y)$ for $x,y\in\M_i$, and hence a Hilbert
$\N$--$\N$ bimodule in the sense of \cite{lance1995hilbert}.
The orthogonal direct sum
$\bigoplus_i \M_i$ is again an $\N$--$\N$ bimodule with inner product
$\langle\oplus_i x_i,\oplus_i y_i\rangle_{\N}=\sum_i\langle x_i,y_i\rangle_{\N}$.
On simple tensors, the interior tensor product
$\mathring{\M}_{i_1}\otimes_\N\cdots\otimes_\N\mathring{\M}_{i_n}$ carries the inner product
\begin{equation}\label{eq:tensor-ip}
\langle x_1\otimes\cdots\otimes x_n,  y_1\otimes\cdots\otimes y_n\rangle_{\N}
=E_{i_n}  \Big(x_n^* E_{i_{n-1}}  \big(\cdots E_{i_1}(x_1^*y_1)\cdots\big) y_n\Big).
\end{equation}
The AFP $\M_1*_{\N}\M_2$ can be realized as the completion of the
$*$-algebra of reduced words.

\begin{equation}\label{eq:reduced-word}
  \N \oplus
  \bigoplus_{n=1}^\infty
    \bigoplus_{\substack{i_1,\dots,i_n\in\{1,2\}\\ i_k\neq i_{k+1}}}
    \mathring{\M}_{i_1}\otimes_\N \cdots \otimes_\N \mathring{\M}_{i_n}.
\end{equation}
This decomposition inherits the structure of a Hilbert $\N$--$\N$ bimodule by discussion around \cref{eq:tensor-ip}.

There are natural embeddings
\begin{equation}\label{eq:pi_i}
\pi_i:\M_i\hookrightarrow \M_1*_{\N}\M_2,\qquad i=1,2.
\end{equation}
The canonical trace on the AFP, denoted by $\wt\tau:=\tau_1*\tau_2$, is defined on reduced words.
By \cite[Theorem~4.2]{takesaki_theory_2003-2}, there exists a unique $\wt\tau$-preserving conditional
expectation $\wt E:\M_1*_{\N}\M_2\to\N$ such that:
\begin{enumerate}[(i)]
  \item $\pi_1|_{\N}=\pi_2|_{\N}=\Id_{\N}$;
  \item $\wt\tau\circ\wt E=\wt\tau$;
  \item $\wt E\circ\pi_i=E_i$, $i=1,2$;
  \item $\wt\tau|_{\N}=\tau_1|_{\N}=\tau_2|_{\N}$;
  \item (\emph{Freeness with amalgamation} \cite{voiculescu_free_1992})
  if $i_1,\dots,i_n\in\{1,2\}$ with $i_k\neq i_{k+1}$ and $m_k\in\mathring{\M}_{i_k}$, then
  \[
  \wt\tau(\pi_{i_1}(m_1)\cdots\pi_{i_n}(m_n))=0,
  \qquad
  \wt E(\pi_{i_1}(m_1)\cdots\pi_{i_n}(m_n))=0;
  \]
  \item $\M_1*_{\N}\M_2$ is generated by $\pi_1(\M_1)\cup\pi_2(\M_2)$.
\end{enumerate}
The above properties give the abstract characterization of AFP.
\paragraph{The AFP Fock space representation}
We also construct the AFP Fock space with non-scalar $\N$.
In the scalar case where $\N=\bC\xi$, this reduces to the free Fock space representation \Cref{eq:free-fock}.
With the notation above, the GNS Hilbert space of the AFP $(\M_1*_{\N}\M_2,\wt\tau)$ is
\[
\H_{12}:=\overline{\M_1*_{\N}\M_2}^{\|\cdot\|_2},
\qquad
\|x\|_2=\wt\tau(x^*x)^{1/2}.
\]
with scalar inner product induced from the $\mathcal{N}$-valued inner product:
\begin{equation}\label{eq:innerproduct-trace}
\langle x,y\rangle := \wt\tau(\wt E(x^*y))
= \wt\tau(\langle x,y\rangle_\N) \in \bC .
\end{equation}
One can then write $\H_{12}$\footnote{This construction is equivalent to introducing Connes fusion \cite{takesaki_theory_2003-2}, but we avoid introducing this machinery here for simplicity.} as \cite{brown_-algebras_2008}:
\begin{equation}\label{eq:hilbert-free}
\H_{12}
= L^{2}(\N)
\oplus
\bigoplus_{n=1}^\infty
\bigoplus_{\substack{i_{1},\dots,i_{n}\in\{1,2\}\\ i_{k}\neq i_{k+1}}}
L^{2}     (\mathring{\M}_{i_{1}}\otimes_{\N}\cdots\otimes_{\N}\mathring{\M}_{i_{n}}   ).
\end{equation}
where $\oplus$ denotes the orthogonal direct sum with respect to the inner product in \eqref{eq:innerproduct-trace}.

Each inclusion $\pi_i:\M_i\hookrightarrow \M_1*_{\N}\M_2\subset B(\H_{12})$, $i=1,2$,
defines a left regular action on simple tensors, which form a dense set, given by \cite[Lemma 4.7.3]{brown_-algebras_2008}
\[
\pi_i(x)(n\xi_\N)=(x-E_i(x))\otimes_\N n\xi_\N+E_i(x)n\xi_\N,\qquad n\in\N,
\]
where $\xi_\N$ denotes the GNS vector for $L^2(\N)=\overline{\N \xi_\N}$. For $\zeta_1\otimes\cdots\otimes\zeta_n \in \mathring{\M}_{i_1}\otimes_{\N}\cdots\otimes_{\N}\mathring{\M}_{i_n}, i_1,\dots,i_n\in\{1,2\}$,
\begin{equation}\label{eq:action}
\pi_i(x)(\zeta_1\otimes\cdots\otimes\zeta_n)
  =
  \begin{cases}
    \bigl(x \zeta_1-E_i(x\zeta_1)\bigr)\otimes\zeta_2\otimes\cdots\otimes\zeta_n
    +E_i(x  \zeta_1)  \zeta_2\otimes\cdots\otimes\zeta_n,
    & i_1=i,\\
    (x-E_i(x))\otimes  \zeta_1\otimes\cdots\otimes\zeta_n+E_i(x)\zeta_1\otimes\cdots\otimes\zeta_n,
    & i_1\neq i,
  \end{cases}
\end{equation}
In words, if the leading index equals $i$, $\pi_i(x)$ acts by left multiplication on the first tensor leg, with the usual centering correction on $\mathring{\M}_i$; otherwise, it creates a new tensor leg on the left up to centering correction\footnote{Under the identification obtained by regrouping the decomposition
\eqref{eq:reduced-word} of
$\M_1 *_\N \M_2$ as $\M_i\otimes_\N \cdots$, the representations coincide with the left regular actions:
$\pi_i(x) = x \otimes \mathbf{1}$.}.
Each $\pi_i$ is a faithful normal $*$-homomorphism with $\pi_i|_{\N}=\Id_{\N}$. Thus \cref{eq:pi_i} extends to the \emph{GNS representation} after completion.
\section{Poincar\'e inequality for trace-symmetric QMS}\label{sec:ppp-tracial}
In this section we prove \Cref{main:tracial}.
The argument proceeds as follows:
We introduce the swap-operator construction as a noncommutative symmetrization,
state auxiliary lemmas for the associated derivation,
establish a generalized Klein inequality giving a convex chain rule,
and finally prove the $\Poin(p,p)$ inequality and its non-self-adjoint extension. To ease notation, in this section, we write p-norms $\|\cdot\|_{L^p(\M,\tau)}$ as $\|\cdot\|_p$.
\subsection{Swap operator construction}\label{sec:swap}
Specifically, we are working in the case where $(\M_1,\tau_1)$ and $(\M_2,\tau_2)$ are two copies of $(\M,\tau)$ containing the von Neumann subalgebra $\mathcal N$.
We use the shorthand $\H_{ij}:=L^2(\M_i *_{\N} \M_j,\wt\tau)$ and define a swap operator
\[
U \colon \H_{12} \to \H_{21}
\]
that exchanges the roles of the two legs $1   \leftrightarrow 2$ in the Hilbert space generated by the AFP $\M_1 *_{\N} \M_2$.
As a prelude, we first display the commuting diagram that organizes the construction.
The swap unitary $U$ implementing $\Ad_U$ is constructed later in this subsection.
\begin{figure}[htbp]
    \centering
\begin{tikzcd}[column sep=5.5em,row sep=4em]
  & \M_{1}
      \arrow[dr,"E_{1}"]
      \arrow[dl,hook',bend left=15,"\pi_{1}"']
      \arrow[dl,bend right=45,<-,"\cE_{1}"']
  \\[2pt]
  \M_{1}\ast_{\mathcal N}\M_{2}
      \arrow[rr,"\exists! \wt E" description]
      \arrow[out=155,in=205,loop,looseness=6,swap,"\Ad_U"]
     & & \mathcal N
     \arrow[out=40,in=-40,loop,looseness=8,swap,"\Ad_U|_{\N}=\Id|_{\N}" {right}]
  \\[2pt]
  & \M_{2}
      \arrow[ur,"E_{2}"']
      \arrow[ul,hook,bend right=15,"\pi_{2}"]
      \arrow[ul,bend left=45,<-,"\cE_{2}"']
\end{tikzcd}
\caption{Roadmap for the AFP construction and the swap automorphism.}
\label{fig:amalgamated-diagram}
\end{figure}

\noindent
Here, each arrow denotes
\begin{enumerate}[(i)]
  \item $\pi_i\colon \M_i\hookrightarrow \M_{1}\ast_{\N}\M_{2}$ is the canonical inclusion.
  \item $\cE_i\colon \M_{1}\ast_{\N}\M_{2}\to \M_i$ is the canonical conditional expectation onto $\M_i$, characterized by
  \[
    \cE_i\circ\pi_i=\Id_{\M_i},
    \qquad
\cE_i(\pi_j(x))=E_j(x) \quad (j\neq i).  \]
  \item $E_i\colon \M_i\to \N$ is the given $\tau_i$-preserving conditional expectation.
  \item The two triangles commute, i.e.\ $E_i\circ\cE_i=\wt E$.
  \item $\wt E\colon \M_{1}\ast_{\N}\M_{2}\to \N$ is the unique conditional expectation whose restriction to $\pi_i(\M)$ equals $E_i$.
  \item $\Ad_U\colon \M_1\ast_{\N}\M_2\to\M_1\ast_{\N}\M_2$ is the spatial $^*$-automorphism implemented by the swap unitary $U$, satisfying
  \[
    \Ad_U\circ\pi_1=\pi_2,
    \qquad
    \Ad_U\circ\pi_2=\pi_1,
    \qquad
    \Ad_U|_{\N}=\Id|_{\N}.
  \]
\end{enumerate}
Expanding the canonical orthogonal decomposition \cref{eq:hilbert-free},
\[
\begin{aligned}
  \H_{12}
  &=
  L^2(\N)
  \oplus
  \bigl( L^2(\mathring{\M}_1) \oplus L^2(\mathring{\M}_2) \bigr) \\
  &\quad \oplus
  \bigl( L^2(\mathring{\M}_1 \otimes_{\N} \mathring{\M}_2)
        \oplus
        L^2(\mathring{\M}_2 \otimes_{\N} \mathring{\M}_1) \bigr) \\
  &\quad \oplus
  \bigl( L^2(\mathring{\M}_1 \otimes_{\N} \mathring{\M}_2 \otimes_{\N} \mathring{\M}_1)
        \oplus \cdots \bigr).
\end{aligned}
\]
Similarly, the decomposition of $\H_{21}$ swaps $1   \leftrightarrow2$.
Both decompositions are canonically isometric. In the symmetric case $\M_1=\M_2=\M$, the swap operator is defined by
\[
\begin{aligned}
  U\colon
  L^2(\mathring{\M}_{i_1} \otimes_{\N} \cdots \otimes_{\N} \mathring{\M}_{i_n})
   &\longrightarrow
   L^2(\mathring{\M}_{3-i_1} \otimes_{\N} \cdots \otimes_{\N} \mathring{\M}_{3-i_n}), \\
  \eta_1 \otimes \cdots \otimes \eta_n
   &\longmapsto
   \zeta_1 \otimes \cdots \otimes \zeta_n,
\end{aligned}
\]
with $\zeta_k=\eta_k$ regarded in the corresponding copy of $\o{\M}_{3-i_k}$.
On the vacuum sector we set $U|_{L^2(\N)}=\Id|_{L^2(\N)}$.
This implies a unitary isomorphism $\H_{12}\simeq \H_{21}$ implementing the label exchange $1   \leftrightarrow 2$.

\smallskip
\begin{proposition}
The swap operator $U$ is unitary.
\end{proposition}

\begin{proof}
It suffices to verify norm preservation on simple tensors.
Let
\[
\eta=\eta_1\otimes\cdots\otimes\eta_n
\in \mathring{\M}_{i_1}\otimes_{\N}\cdots\otimes_{\N}\mathring{\M}_{i_n},
\]
and set $\zeta_k=\eta_k\in\mathring{\M}_{3-i_k}$.
Using the $L^2$ inner product \cref{eq:innerproduct-trace} and $\wt\tau|_{\N}=\tau|_{\N}$,
\begin{align*}
\|U\eta\|_2^2
&=\langle U\eta,U\eta\rangle
=\wt\tau(\wt E((U\eta)^*(U\eta))) \\
&=\tau(E_{3-i_n}(\zeta_n^*
E_{3-i_{n-1}}(\cdots E_{3-i_1}(\zeta_1^*\zeta_1)\cdots)\zeta_n)) \\
&=\tau(E_{i_n}(\eta_n^*
E_{i_{n-1}}(\cdots)\eta_n))
=\|\eta\|_2^2,
\end{align*}
where the last equality follows by relabeling $3-i_k\mapsto i_k$.
Thus $U$ extends to an isometry $\H_{12}\to\H_{21}$.
The inverse construction gives an isometry $\H_{21}\to\H_{12}$, hence $U$ is unitary.
\end{proof}

\begin{proposition}\label{thm:swap}
Let $U$ be the swap unitary on $\H_{12}$.
Then $\Ad_U(T)=UTU^*$ restricts to a $*$-automorphism of $\M_1*_{\N}\M_2$ satisfying:
\begin{enumerate}[(i)]
\item $\Ad_U(\pi_i(x))=\pi_{3-i}(x)$ for all $x\in\M$.
\item $\Ad_U$ commutes with $\wt E$.
\item $\Ad_U$ preserves $\wt\tau$.
\end{enumerate}
\end{proposition}

\begin{proof}
Since $U$ is unitary, $\Ad_U$ defines a spatial $*$-isomorphism between
$B(\H_{12})$ and $B(\H_{21})$ \cite[Section 3]{Takesaki_1979}.
Identifying $\H_{21}$ with $\H_{12}$ via $U$, we view $\Ad_U$ as an automorphism
of $B(\H_{12})$.

\emph{(i)}
It suffices to check the action on simple tensors.
Let $\eta=\eta_1\otimes\cdots\otimes\eta_n$ with $i_1=i$.
Using \cref{eq:action} and $U|_{L^2(\N)}=\Id$,
\begin{align*}
U\pi_i(x)\eta
&=(x\zeta_1-E_{3-i}(x\zeta_1))\otimes\cdots\otimes\zeta_n
+E_{3-i}(x\zeta_1)\zeta_2\otimes\cdots\otimes\zeta_n \\
&=\pi_{3-i}(x)(\zeta_1\otimes\cdots\otimes\zeta_n)
=\pi_{3-i}(x)U\eta,
\end{align*}
where $\zeta_k=\eta_k$ viewed in $\mathring{\M}_{3-i_k}$.
Thus $\Ad_U(\pi_i(x))=\pi_{3-i}(x)$.
Since $\M_1*_{\N}\M_2$ is generated by $\pi_1(\M)$ and $\pi_2(\M)$,
$\Ad_U$ restricts to a $*$-automorphism.

\emph{(ii)} Following \Cref{fig:amalgamated-diagram}, \[ \Ad_U  \circ \wt E\circ\pi_i=\wt E \circ \pi_i=\wt E \circ \pi_{3-i}=\wt E \circ \Ad_U \circ \pi_i, \] since $U|_{L^2(\N)}=\Id|_{L^2(\N)}$ and $\M_1=\M_2=\M$. As $\M_1*_\N \M_2$ is generated by $\pi_1(\M)$ and $\pi_2(\M)$ and both $\Ad_U$ and $\wt E$ are normal, the relation extends to all of $\M_1*_\N \M_2$. Thus $\Ad_U\circ \wt E=\wt E =\wt E\circ \Ad_U$.
\emph{(iii)} Because $\wt\tau=\wt\tau\circ\wt E$, \[ \wt\tau(\Ad_U(T)) = \wt\tau(\wt E(\Ad_U(T))) = \wt\tau(\wt E(T)) = \wt\tau(T), \qquad T\in\M_1*_\N\M_2. \] Hence $\Ad_U$ swaps the two embeddings, commutes with $\wt E$, and preserves $\wt\tau$.
\end{proof}

We will see in \Cref{sec:ppp-proof} that the swap unitary $U$
is the essential ingredient for the symmetrization technique,
replacing the product measure in the commutative case
\cite[Section~4.1]{huang2021poincare}.

\subsection{Auxiliary lemmas}
We first show several auxiliary lemmas used in the proof of \Cref{main:tracial}.
They state the relevant estimates in terms of a noncommutative derivation,
replacing the product-measure construction of doubling-up variables techniques from \cite{huang2021poincare} in the semi-commutative case.
We define the noncommutative derivation
\[
  \delta(x):=\pi_1(x)-\pi_2(x)\in\M_1*_\N\M_2 .
\]
\begin{lemma}\label{lem:contraction}
Let $x\in\M$.
Then for any $1\le p<\infty$,
\[
  \|x-E_1(x)\|_p^p \le \|\delta(x)\|_p^p .
\]
\end{lemma}

\begin{proof}
Since the conditional expectation $\cE_1:L_p(\M_1*_\N\M_2)\to L_p(\M_1)
$ is contractive \cite{takesaki_theory_2003-2}, we have
\[
\|\cE_1(\delta(x)-E(\delta(x)))\|_p\le \|\delta(x)-E(\delta(x))\|_p .
\]
Moreover,
\[
\cE_1(\delta(x))
=\cE_1(\pi_1(x))-\cE_1(\pi_2(x))
=x-E_1(x),
\]
and since $E\circ\pi_i=E_i$ with $E_1=E_2$, one has $E(\delta(x))=0$.
The claim follows.
\end{proof}
The following lemma gives the noncommutative version of the symmetrization technique:
\begin{lemma}\label{lem:odd}
Let $x\in\M_{\mathrm{sa}}$, and let $\delta(x)$ be as above and $U$ the swap unitary satisfying
$\Ad_U\circ\pi_1=\pi_2$ and $\Ad_U\circ\pi_2=\pi_1$,
let $h:\R\to\R$ be an odd function. Then
\begin{enumerate}[(i)]
  \item $\Ad_U(h(\delta(x)))=-h(\delta(x))$,
  \item $\wt\tau(h(\delta(x)))=0$,
  \item $E(h(\delta(x)))=0$,
\end{enumerate}
where $\wt\tau$ is the canonical trace on $\M_1*_\N\M_2$.
\end{lemma}

\begin{proof}
By \Cref{thm:swap}, $\Ad_U(\delta(x))=-\delta(x)$.
Since $h$ is odd, $h(-\delta(x))=-h(\delta(x))$.
The $U$-invariance of $\wt\tau$ and $E$ then implies the stated identities.
\end{proof}

\begin{lemma}\label{lem:factor-2}
Let $\wt L$ and $\wt \Gamma$ be the Lindblad generator and gradient form on
$\M_1*_\N\M_2$ induced from $\M$ as in \cite{boca_completely_1993}, characterized by
\[\wt L\bigl(\pi_i(x)\bigr)
= \pi_i\bigl(L_i(x)\bigr),
\qquad
\wt\Gamma\bigl(\pi_i(x),\pi_i(x)\bigr)
= \pi_i\bigl(\Gamma_i(x,x)\bigr),
\quad i=1,2,\]
where $L_i=L$ and $\Gamma_i=\Gamma$ are the corresponding operators on $\M$.
Then for every $x\in\M$ and $p\ge2$,
\[
  \|\wt\Gamma(\delta(x),\delta(x))\|_{p/2}
  \le 2\|\Gamma(x,x)\|_{p/2}.
\]
\end{lemma}

\begin{proof}
A direct computation shows
\[
\wt \Gamma(\delta(x),\delta(x))
=\pi_1(\Gamma(x,x))+\pi_2(\Gamma(x,x)).
\]
Since $\pi_i$ preserves the $L^{p/2}$-norm and the trace,
the inequality follows immediately.
\end{proof}

\subsection{Generalized Klein inequalities}\label{sec:Klein}
The key analytic input for the tracial $\Poin(p,p)$ inequality is a convexity estimate
relating functional calculus and the gradient form.
This is provided by Markov dilation (\Cref{app:chain}) and a generalized Klein inequality, which is used to obtain a convex-chain rule
for the Dirichlet form.
We recall the tracial version below and derive the corresponding gradient estimate.
\medskip
\begin{lemma}[Generalized Klein inequality for tracial von Neumann algebra]
\label{lem:klein}
Let $(\M,\tau)$ be a von Neumann algebra as in the last section.  Let
$
x,y\in\M_{\mathrm{sa}}
$
be bounded self-adjoint elements, and let $\varphi:\R\to\R$ be a continuously differentiable function whose squared derivative
$
\psi = (\varphi')^2
$
is convex on $\R$.  Then
\[
    \tau   \bigl[(\varphi(x)-\varphi(y))^2\bigr]
 \le
\frac12
 \tau   \Bigl[(x-y)^2\bigl(\psi(x)+\psi(y)\bigr)\Bigr].
\]
\end{lemma}
\begin{proof}
By \cite[Prop.~3]{petz_survey_1994}, the claim holds in finite dimensions, and the general bounded case follows by spectral approximation and norm continuity.
\end{proof}
\begin{proposition}[Convex-chain estimate for gradient form]
\label{prop:convex-chain}
Assume the hypotheses of Section~\ref{sec:preliminaries} and lemma \ref{lem:klein},
for every self-adjoint element
$x\in\Dom(L)\subset L_{2}(\M,\tau)$,
\[
  \mathscr E\bigl(\varphi(x)\bigr)
   \le
  \tau   \bigl(\Gamma(x,x)  \psi(x)\bigr),
\]
where
$\mathscr E(f):=\tau   \bigl(\Gamma(f,f)\bigr)$ is the Dirichlet form
associated with $L$.
\end{proposition}
We defer the proof to \Cref{app:chain}, where the main technique is Markov dilation that generalizes the Markov processes in \cite{huang2021poincare}.
We use the main dilation theorem of
\cite{JungeRicardShlyakhtenko2025}, which states that every weak$^*$-continuous
semigroup of normal, unital, completely positive maps that is self-adjoint on
$L^2(\M,\tau)$ admits a Markov dilation in the sense of
Definition~\ref{def:dilation}.
\subsection{Proof of trace-symmetric \texorpdfstring{$\Poin(p,p)$}{Poin(p,p)}}\label{sec:ppp-proof}
We now prove \Cref{main:tracial}. We begin by recalling its statement.
\begin{theorem*}
Let $(\M,\tau)$ be a tracial von Neumann algebra and $(T_t)_{t\ge0}=e^{-tL}$ a $\tau$--symmetric quantum Markov semigroup with fixed-point algebra $\N$ and conditional expectation $E:\M\to\N$.
If $\{T_t\}_{t\ge0}$ has a spectral gap $\alpha>0$, then for every self-adjoint $x$ in the domain of $L$ and $p=2$ or $p\ge3$, it satisfies the Poincar\'e inequality $\Poin(p,p)$ with constant $\tfrac{p}{\sqrt{2\alpha}}$, i.e.,
\[
\|x-E(x)\|_p  \le  \frac{p}{\sqrt{2\alpha}} \|\Gamma(x,x)^{1/2}\|_p.\tag{$\Poin(p,p)$}
\]
\end{theorem*}
\begin{proof}
For $x\in\M_{sa}$, \Cref{lem:contraction} gives
\begin{equation}\label{eq:proof1}
\|x-E_1(x)\|_p^p \le \|\delta(x)\|_p^p,
\qquad
\delta(x):=\pi_1(x)-\pi_2(x)\in(\M_1*_\N\M_2)_{sa}.
\end{equation}

Define the odd function $\varphi(s):=\operatorname{sgn}(s)|s|^{p/2}$ (cf.\ \cite{huang2021poincare}).
Its squared derivative is
\[
\psi(s):=(\varphi'(s))^2=\Big(\frac p2\Big)^2|s|^{p-2},
\]
which is convex for $p=2$ and $p\ge3$.
Since $\varphi$ is odd, \Cref{lem:odd} implies
$E(\varphi(\delta(x)))=0$, and hence
\[
\|\varphi(\delta(x))-E(\varphi(\delta(x)))\|_2^2
=\wt\tau(\varphi(\delta(x))^*\varphi(\delta(x)))
=\|\delta(x)\|_p^p.
\]

With the spectral-gap assumption \eqref{eq:gap}
(extended to $\M_1*_\N\M_2$ by \Cref{cor:gap}),
\[
\alpha\|\delta(x)\|_p^p
\le
\|\wt\Gamma(\varphi(\delta(x)),\varphi(\delta(x)))^{1/2}\|_2^2
=\wt\tau(\varphi(\delta(x))^*\wt L(\varphi(\delta(x))))
=:\mathscr E(\varphi(\delta(x))).
\]
By \Cref{prop:convex-chain} (convex chain rule),
\[
\mathscr E(\varphi(\delta(x)))
\le
\wt\tau(\wt\Gamma(\delta(x),\delta(x))\psi(\delta(x)))
=\Big(\frac p2\Big)^2
\wt\tau(\wt\Gamma(\delta(x),\delta(x))|\delta(x)|^{p-2}).
\]

Applying Hölder's inequality,
\[
\alpha\|\delta(x)\|_p^p
\le
\Big(\frac p2\Big)^2
\wt\tau(\wt\Gamma(\delta(x),\delta(x))^{p/2})^{2/p}
\|\delta(x)\|_p^{p-2}.
\]
Dividing by $\|\delta(x)\|_p^{p-2}$ and applying \Cref{lem:factor-2} gives
\[
\|\delta(x)\|_p^2
\le
\frac{p^2}{4\alpha}
\wt\tau(\wt\Gamma(\delta(x),\delta(x))^{p/2})^{2/p}\le
\frac{p^2}{2\alpha}
\tau(\Gamma(x,x)^{p/2})^{2/p}.
\]
Taking square roots and using \eqref{eq:proof1} completes the proof:
\[
\|x-E(x)\|_p
\le
\frac p{\sqrt{2\alpha}}
\|\Gamma(x,x)^{1/2}\|_p.
\]
\end{proof}
\begin{remark}[Missing Parameters \cite{huang2021poincare}] \Cref{main:tracial} also holds for $p\in(2,3)$ with an extra factor of $\sqrt{2}$ on the right hand side through applying a variant of \Cref{lem:klein} that only needs $\psi$ to be monotone.

\end{remark}
We can also generalize to non-self-adjoint elements in $\M$ by the following corollary.
We note that in the non-tracial setting, the one-sided estimate in \Cref{eq:non-self-adjoint} may be useful in the analysis of quantum algorithms.
\begin{cor}[Non-self-adjoint $\Poin(p,p)$]\label{cor:non-self-adjoint-ppp}
Under the assumptions of \Cref{main:tracial},
for every $p=2$ and $p\geq 3$ and $x\in\M$:
\[
\|x-E(x)\|_{p} \leq \frac{p}{\sqrt{2\alpha}}  \bigl(\|\Gamma(x,x)\|^{1/2}_{p/2}+\|\Gamma(x^*,x^*)\|^{1/2}_{p/2}\bigr).
\]
\end{cor}
\begin{proof}
Since $x=\Re x+\mathrm i\Im x$ with $\Re x,\Im x\in\M_{sa}$, the triangle inequality in $L_p$ gives
\[
\|x-E(x)\|_p
\le
\|\Re x-E(\Re x)\|_p
+
\|\Im x-E(\Im x)\|_p.
\]
We apply \Cref{main:tracial} to $\Re x$ and $\Im x$,
\[
\|\Re x-E(\Re x)\|_p
+
\|\Im x-E(\Im x)\|_p
\le
\frac p{\sqrt{2\alpha}}
\bigl(
\|\Gamma(\Re x,\Re x)\|_{p/2}^{1/2}
+
\|\Gamma(\Im x,\Im x)\|_{p/2}^{1/2}
\bigr).
\]
Since
\[
\|u\|_{\Gamma,p}:=\|\Gamma(u,u)\|_{p/2}^{1/2}
\]
is a seminorm and
\[
\Re x=\frac{x+x^*}{2},
\qquad
\Im x=\frac{x-x^*}{2\mathrm i},
\]
its triangle inequality gives
\[
\|\Gamma(\Re x,\Re x)\|_{p/2}^{1/2}
+
\|\Gamma(\Im x,\Im x)\|_{p/2}^{1/2}
\le
\|\Gamma(x,x)\|_{p/2}^{1/2}
+
\|\Gamma(x^*,x^*)\|_{p/2}^{1/2}.
\]
Combining the above estimates gives the stated corollary.
\end{proof}

\section{Poincar\'e inequality for GNS-detailed balance QMS}\label{sec:GNS}
Previously, we studied trace-symmetric quantum Markov semigroups and proved the
$\Poin(p,p)$ inequality (\cref{main:tracial}) under the assumptions of a spectral gap
on a tracial von Neumann algebra $(\M,\tau)$.
In this section, we extend this result to von Neumann algebras that are not equipped
with a trace, and to QMSs satisfying the corresponding GNS detailed-balance
condition with respect to a non-tracial state.
The main tools are Haagerup's reduction method
\cite{haagerup2009reductionmethodnoncommutativelpspaces}
and Kosaki's $L^p$ interpolation theorem
\cite{kosaki_applications_1984}.
We also formulate the results using Haagerup's functorial $L^p$ construction. To set notation, we first recall results from Haagerup's $L^p$ spaces and Kosaki's $L^p$ interpolation spaces.
\subsection{Haagerup's and Kosaki's \texorpdfstring{$L^p$}{Lp} spaces}\label{sec:lp-spaces}
\subsubsection{Haagerup's \texorpdfstring{$L^p$}{Lp} spaces}\label{sec:H-Lp}
We give a concise introduction to Haagerup's $L^p$ spaces \cite[Section~1.2]{haagerup2009reductionmethodnoncommutativelpspaces}\cite{haagerup1979lpspaces} for general von Neumann algebras.
Let $\M$ be a von Neumann algebra and $\phi$ a faithful normal semifinite weight on $\M$.  Denote by $ \sigma_t^\phi=\sigma_t\ (t\in\R)$
its modular automorphism group, and form the crossed product \cite[Section~X]{takesaki_theory_2003-2},
$
  \cR=\M\rtimes_{\sigma}\R
   = \bigl\{\pi_\sigma(\M)\cup\lambda(\R)\bigr\}''
$
, where $\pi_\sigma$ and $\lambda$ act on $L^2(\M,\phi)\otimes L^2(\R,dt)$.
It is well known that $\cR$ is semifinite and admits a unique normal semifinite faithful trace $\tau$ (the \emph{canonical trace}) characterized by the scaling
\begin{equation}\label{eq:canonical-trace}
  \tau\circ\wh\sigma_t = e^{-t} \tau,
  \quad t\in\R,
\end{equation}
where $\wh\sigma$ is the dual action on $\cR$.

We define and summarize some useful properties of Haagerup's $L^p$--spaces \cite{haagerup1979lpspaces}.
\begin{definition}[Haagerup's $L^p$ Spaces]\label{def:Hlp}
    With the notation above, for $0 < p \le \infty$, let $\wt\cR$ be the space of $\tau$-measurable\footnote{An operator $a$ on $L^2(\M,\phi)\otimes L^2(\R,dt)$ affiliated with $\cR$ is called $\tau$-measurable \cite{Terp1981} if it is measurable with respect to $(\cR,\tau)$.} operators affiliated with $\cR$. For $0 < p \le \infty$, the Haagerup $L^p$-space is defined by
\[
  L^p(\M)
  :=\bigl\{ x\in \wt \cR \big| \wh\sigma_t(x)=e^{-t/p} x,\quad t\in\R\bigr\}.
\]
\end{definition}
\begin{proposition} \label{prop:h-lp}
From \Cref{def:Hlp} of the Haagerup $L^p(\M)$,
we have the following properties:
\begin{enumerate}[label=(\roman*)]
  \item $L^p(\M)$ is a closed, self--adjoint $\M$--bimodule in $\wt \cR$, stable under $x\mapsto x^*$ and $x\mapsto|x|$.
  \item
$L^p(\M)$ is linearly spanned by the positive part  $
  L^p(\M)_+
  := L^p(\M)\cap\wt\cR_+
$.

  \item $L^\infty(\M)=\M$, since the fixed-point algebra of $\wh \sigma$, $\cR^{\wh\sigma}=\pi_\sigma(\M)\cong\M$.
  \item The map $\M_*^+ \ni \omega\mapsto h_\omega\in L^1(\M)_+$ is an isometric isomorphism, and the dual weight $\wh \omega$ on $\cR$ is equal to $\tau(h_\omega \cdot).$
  \item By linear extension of the above bijection, we have an isometric isomorphism $\mathrm{tr}: \M_* \ni \omega\mapsto h_\omega\in L^1(\M)$, by $\mathrm{tr}(h_\omega):=\omega(1)$. Then we have $$\|{h_\omega}\|_{1}:=\mathrm{tr}(|h_\omega|)=\mathrm{tr}(h_{|\omega|})=|\omega|(1)=\|{\omega}\|.$$
  \item Let $x\in\wt\cR$ with polar decomposition $x=u|x|$, for every $p\in[1,\infty)$, $x\in L^p(\M)$  if and only if $|x|^p\in L^1(\M)$ and $u\in\M$. Define
\[
      \|x\|_{L^p(\M)}
  := \bigl\| |x|^p\bigr\|_1^{1/p}
  \quad(0<p<\infty),
  \qquad
  \|x\|_\infty := \|x\|.
\]

 \item The multiplication in $\wt\cR$, restricts to a bounded bilinear map
\[
  L^p(\M)\times L^{s}(\M) \longrightarrow L^r(\M),
  \quad  \tfrac1p + \tfrac1{s}=\tfrac1r,
\]
is of norm 1. Equivalently, Hölder's inequality extends to Haagerup's $L^p$ spaces.  In particular, if $1/p+1/q=1$, then the pairing
\[
  L^p(\M)\times L^{q}(\M) \ni (x,y) \mapsto \mathrm{tr}(xy)\in \mathbb{C}
\]
identifies $L^{q}(\M)$ isometrically with the dual of $L^p(\M)$ (if $p\neq\infty$), and one has the tracial property
\[
  \mathrm{tr}(xy)=\mathrm{tr}(yx),
  \quad x\in L^p(\M), y\in L^{q}(\M).
\]
\end{enumerate}
\end{proposition}

\begin{remark}
$L^p(\M)$ is independent of the choice of the weight $\phi$ \cite{haagerup1979lpspaces,Terp1981}(up to canonical isometric isomorphism), so it is customary to write simply $L^p(\M)$.  When $\phi=\tau_0$ is a trace, one recovers the Dixmier--Segal noncommutative $L^p$ spaces $L^p(\M,\tau_0)$.
\end{remark}
\subsubsection{Kosaki's \texorpdfstring{$L^p$}{Lp} spaces} \label{sec:K-Lp}
In this section, we give a brief introduction to Kosaki's $L^p$ spaces based on \cite[Appendix~A.8]{Hiai_2021}, where the noncommutative $L^p$ spaces are defined as complex interpolation spaces between $\M$ and $L^1(\M)$. It also turns out that there is a natural isomorphism \cite[Theorem~9.1]{kosaki_applications_1984} of Kosaki's and Haagerup's $L^p$ spaces. Using the same notation as the previous paragraph with an extra constraint of $\M$ \emph{$\sigma$-finite} and $\phi$ be a faithful state so that $\phi(\cdot)=\mathrm{tr}(h_\phi \cdot)$, the first step towards interpolation is to construct a pair of two compatible Banach spaces from $\M$ and $\M_*$. One choice is the injective maps
\[
    \iota_\eta^\phi: \M\ni a\mapsto h_\phi^{\eta} a h_\phi^{1-\eta}\in L^1(\M)\quad\text{for}\quad 0\le\eta\le 1.
\]
We drop the subscript $\phi$ whenever the state is clear from context.
Define the norm $\|\iota_\eta(a)\|:=\|a\|$ on $\iota_\eta(\M)\subset L^1(\M)$, so that $\M\cong\iota_\eta(\M)$. Therefore, the pair $(\iota_\eta(\M),L^1(\M))$ form a pair of compatible Banach spaces.
\begin{definition}[Kosaki's $L^p$ spaces {\cite[Def.~A.60]{Hiai_2021}}]\label{def:k-lp}
With the constructions as above,
let $1<p<\infty$ and $0\le\eta\le1$.  Kosaki's $L^p$ space $L^p(\M,\phi)_\eta$ with respect to $\phi$ is defined to be the complex interpolation space \cite[Section~I.1]{kosaki_applications_1984}
\[  C_{1/p}\bigl(\iota_\eta(\M), L^1(\M)\bigr),
\]
equipped with the complex interpolation norm
$\|\cdot\|_{p,\phi,\eta} (=\|\cdot\|_{C_{1/p}})$, Moreover, one sets
\[
  L^1(\M,\phi)_\eta := L^1(\M),
  \quad
  L^\infty(\M,\phi)_\eta := \iota_\eta(\M).
\]
\end{definition}
In particular, $  L^p(\M,\phi)_0,
  L^p(\M,\phi)_1,
  L^p(\M,\phi)_{1/2}$
are called the \emph{left}, \emph{right} and \emph{symmetric} $L^p$--spaces, respectively. We also establish the isomorphism of different $L^p$ spaces.
\begin{theorem}[Isomorphism of Kosaki's and Haagerup's $L^p$ spaces ]\label{thm:iso-h-k}
Let $1\le p\le\infty$ and let $q$ be its conjugate exponent, $1/p+1/q=1$.  Then for each $0\le\eta\le1$ one has the identification
\[
L^p\bigl(\M,\phi\bigr)_\eta
 =
h_\phi^{\eta/q} L^p(\M) h_\phi^{(1-\eta)/q}
 \subset L^1(\M),
\]
and the map
\begin{equation}\label{eq:h-k-iso}
\iota_{\eta,p} :L^p(\M) \longrightarrow L^p(\M,\phi)_\eta,
\qquad
a\mapsto h_\phi^{\eta/q} a h_\phi^{(1-\eta)/q},
\end{equation}
is an isometry.  Equivalently, for all $a\in L^p(\M)$,
\[
\big\| h_\phi^{\eta/q} a h_\phi^{(1-\eta)/q}\big\|_{p,\phi,\eta}
=\|a\|_{L^p(\M)}.
\]
In particular, $L^p(\M)\cong L^p(\M,\phi)_\eta$ via this isometry.
On the other hand, $a=\iota_{\eta,q}(x)=h_\phi^{\eta/p} x h_\phi^{(1-\eta)/p}$ is dense in $L_{p}(\M)$ for $x\in \M$, so we have
\[
\big\| \iota_{\eta,p} (a)\big\|_{p,\phi,\eta}=\big\| \iota_{\eta}(x)\big\|_{p,\phi,\eta}
=\|a\|_{L^p(\M)}=\|h_\phi^{\eta/p} x h_\phi^{(1-\eta)/p}\|_{L^p(\M)}.
\]
\end{theorem}
In particular, when $\eta=0,1,\tfrac12$ we recover Kosaki's left, right, and symmetric spaces:
\[
  L^p\bigl(\M,\phi\bigr)_L
    :=L^p(\M) h_\phi^{1/q},
  \qquad
  L^p\bigl(\M,\phi\bigr)_R
    :=h_\phi^{1/q} L^p(\M),
\]
\begin{equation}\label{eq:sym-lp}
      L^p\bigl(\M,\phi\bigr)
    :=L^p\bigl(\M,\phi\bigr)_{1/2}
    =h_\phi^{1/2q} L^p(\M) h_\phi^{1/2q},
\end{equation}
with norms
\[
  \bigl\|a h_\phi^{1/q}\bigr\|_{p,\phi,0}
  = \bigl\| h_\phi^{1/q}a\bigr\|_{p,\phi,1}
  = \bigl\| h_\phi^{1/2q} a h_\phi^{1/2q}\bigr\|_{p,\phi,1/2}
  = \|a\|_{L^p(\M)},
  \qquad a\in L^p(\M).
\]
With the above identifications in mind, we also summarize the norms in \Cref{tab:norm}.
\begin{table}[htbp]
\centering
\begin{tabular}{l c}
\hline
\textbf{Spaces} & \textbf{Norms} \\
\hline
Haagerup's $L^p(\M)$ & $\|a\|_{L^p(\M)} = \| |a|^p \|_{1}^{1/p}$ \\[4pt]
Kosaki's space $L^p(\M,\phi)_\eta$ &
$\big\| \iota_{\eta}(x)\big\|_{p,\phi,\eta}
=\|h_\phi^{\eta/p} x h_\phi^{(1-\eta)/p}\|_{L^p(\M)}.$ \\[4pt]
Left $L^p(\M,\phi)_L$ &
$\|x h_\phi\|_{p,\phi,0} = \|xh_\phi^{1/p}\|_{L^p(\M)}$ \\[4pt]

Right $L^p(\M,\phi)_R$ &
$\|h_\phi x\|_{p,\phi,1} = \|h_\phi^{1/p} x\|_{L^p(\M)}$ \\[4pt]

Symmetric $L^p(\M,\phi)_{1/2}$ &
$\|h_\phi^{1/2} x h_\phi^{1/2}\|_{p,\phi,1/2} = \|h_\phi^{\frac{1}{2p}} x h_\phi^{\frac{1}{2p}}\|_{L^p(\M)}$ \\[4pt]
\hline
\end{tabular}
\caption{Norms in Haagerup's and Kosaki's $L^p$ spaces for $a\in L^p(\M)$ and $x\in \M$.}
\label{tab:norm}
\end{table}

\subsection{GNS and KMS detailed balance with respect to a state}\label{sec:GNSsymmetric}
Throughout this section, we work in the GNS representation of $\M$ associated with a normal,
faithful state $\phi$, and use the notation for the modular operator $\Delta$,
modular conjugation $J$, and modular automorphism group $(\sigma_t^\phi)_{t\in\R}$
introduced in \Cref{sec:von-neumann-algebras}.
There are two symmetry conditions \cite[Definition~4.1]{Gao_Junge_LaRacuente_Li_2025}\label{defi:GNS} for QMS relative to the state $\phi$.
\medskip

\begin{definition}[GNS and KMS detailed balance]\label{def:GNS}
A quantum Markov semigroup $(T_t)_{t\geq0}:\M \to \M$ is Gelfand-Naimark-Segal detailed balance with respect to $\phi$ (GNS-$\phi$-detailed balance) if \[\phi(T_t(x)^*y)=\phi(x^*T_t(y)), \quad   \pl x,y\in \M\pl,t\geq 0; \]
We say $T_t$ is Kubo--Martin--Schwinger detailed balance with respect to $\phi$ (KMS-$\phi$-detailed balance) if
\begin{equation}\label{eq:kms}\lan \Delta^{\frac{1}{4}} x\xi_\phi, \Delta^{\frac{1}{4}}T_t(y)\xi_\phi \ran=\lan \Delta^{\frac{1}{4}} T_t(x)\xi_\phi, \Delta^{\frac{1}{4}}y\xi_\phi \ran, \pl   \pl x,y\in \M\pl,t\geq 0. \end{equation}
Correspondingly, we call the pre-adjoint $(T_{t})_*:\M_*\to \M_*$ a GNS- or KMS-$\phi$-detailed balance quantum channel.
\end{definition}

To fix ideas, we restrict attention to the case where $\mathcal{M}$ has a normal, faithful, finite trace $\tau$. The general von Neumann-algebra case follows after the necessary adaptations of notation (see \Cref{sec:H-Lp}).
In the tracial case, every normal state $\phi$ admits a density
$D_\phi\in L^{1}(\M)$ such that $\phi(x)=\tau  \bigl(D_\phi x\bigr)$ for $x\in M$.

With this choice, setting $x\in \M$, the modular operator can be written as
\[\Delta(x)=D_\phi x D_\phi^{-1},\]
and for $s\in\R$ the modular automorphism group \cref{eq:modular-group} takes the form
\[
  \sigma^{\phi}_{s}(x)=D_\phi^{is} x D_\phi^{-is}.
\]
Let $(T_t)_*:L^{1}(\M)\to L^{1}(\M)$ denote the pre-adjoint of the QMS
$T_t$ under the trace duality.  The KMS-$\phi$ detailed balance condition
from \cref{eq:kms} is equivalent to
\begin{equation}\label{eq:kms-2}
  (T_t)_*    \bigl(D_\phi^{1/2}xD_\phi^{1/2}\bigr)
   =
  D_\phi^{1/2} T_t(x) D_\phi^{1/2},
  \qquad x\in\M, t\ge0.
\end{equation}

For $1\le p<\infty$ the symmetric Kosaki $L^{p}$ space (\Cref{sec:K-Lp}) weighted by $\phi$
is obtained by completing $\M$ with respect to the norm \cite{kosaki_applications_1984}, also recalled in \Cref{tab:norm}.
\[
  \|\iota_{1/2}(x)\|_{p,\phi,\frac{1}2{}}
   =
  \bigl\|D_\phi^{1/(2p)} x D_\phi^{1/(2p)}\bigr\|_{p}, \text{ for }
  \|\cdot\|_{p}=\tau(|\cdot|^{p})^{1/p}.
\]
When $p=2$, $L^{2}(\M,\phi)$ is a Hilbert space with KMS-inner product
$
  \| \iota_{1/2}(x)\|^2_{2,\phi,\frac 12}
  =\langle\Delta^{1/4}x\xi_\phi,\Delta^{1/4}x\xi_\phi\rangle.
$
\Cref{eq:kms} shows that $(T_t)_*$ is contractive on $L^{1}(\M,\phi)$. By complex interpolation \cite{kosaki_applications_1984}, it follows that it acts as a contraction on
$L^{p}(\M,\phi)$ (\cref{eq:sym-lp}) for every $1\le p<\infty$.
We also recall from \cite[Proposition~4.2]{Gao_Junge_LaRacuente_Li_2025} some properties of GNS-detailed-balance QMS.

\begin{lemma}\label{lem:cd2}
Let $(T_t)_{t\geq0}:\M\to\M$ be a GNS-$\phi$-detailed-balance quantum Markov semigroup for a normal faithful state $\phi$. Denote by $\N$ its fixed-point algebra (more generally, the multiplicative domain). Then
\begin{enumerate}[(i)]
  \item $\N$ is invariant under $\sigma_t^\phi$.  Hence there exists a $\phi$-preserving normal conditional expectation $E:\M\to\N$.
  \item $T_t|_{\N}$ is an involutive *-automorphism satisfying
    \[
      T_t^2\circ E  =  E\circ T_t^2  =  E,
      \quad
      E\circ T_t  =  T_t\circ E.
    \]
    Moreover, $T_t^2$ is an $\N$-bimodule map,
    $ T_t^2(a x b)=a T_t^2(x) b$  for all $a,b\in\N, x\in\M$.
  \item $T_t$ is isometric on $L_2(\N,\phi)$.  If in addition
    \[
      \bigl\|T_t(\mathrm{id}-E)\colon L_2(\M,\phi)\to L_2(\M,\phi)\bigr\|<1,
    \]
    then
    $
      E  = \displaystyle \lim_{n\to\infty} T_t^{2n} \in B\bigl(L_2(\M,\phi)\bigr).
    $
\end{enumerate}
\end{lemma}
\begin{remark}[Comparisons of GNS and KMS.]\label{rmk:gns-kms}
    Both GNS-- and KMS--detailed balance conditions generalize the classical notion of detailed balance to QMSs. In fact, the GNS--$\phi$--detailed balance condition is strictly \emph{stronger} than the KMS--$\phi$--detailed balance condition: as shown in \cite{Giorgetti:2021wav,Wirth_2024}, a QMS satisfies GNS--$\phi$--detailed balance if and only if \emph{(i)} it satisfies KMS--$\phi$--detailed balance and \emph{(ii)} commutes with the modular automorphism group $\sigma^\phi_t$.
\end{remark}

We will assume the GNS--$\phi$ detailed balance to extend the $\mathrm{PI}(p,p)$ (Theorem~\ref{main:tracial}) beyond the tracial setting to a general non--tracial pair $(\M,\phi)$.
The reason for choosing this stronger condition is its compatibility with Haagerup's reduction in \Cref{sec:reduction}, which allows us to reduce to tracial slices in a larger algebra.
In the sequel, we shall work exclusively under that assumption.

\subsection{Lindbladians, conditional expectations and gradient forms in \texorpdfstring{$L^p$}{Lp} spaces}\label{sec:lind-lp}
We reviewed Haagerup's $L^p(\M)$ and Kosaki's $L^p(\M,\phi)_\eta$ spaces in \Cref{sec:lp-spaces}, assuming that $\M$ is a $\sigma$-finite von Neumann algebra equipped with a normal faithful state $\phi$. As established in \Cref{thm:iso-h-k}, these $L^p$ spaces are naturally isomorphic. In this section, we define $L^p$ versions of Lindbladians, conditional expectation, and associated gradient forms. We will see later that compatibility with the tracial setting holds if and only if the Lindbladian satisfies the GNS-$\phi$-detailed balance condition, which is the necessary regularity assumption in this framework.

\begin{definition}[Lindbladian Generator in Haagerup's $L^p$ Spaces]\label{def:lind-h-lp}
Given the generator $L:\M\to\M$ defined in \cref{eq:lindbladian} on $(\M,\phi)$ and in addition satisfy $\phi$-GNS-detailed balance condition, one can define the $L^p(\M)$ version
\[
    \L_p^H: \Dom(\L_p^H) \subset L^p(\M)\to L^p(\M)
\]
given by dense elements $h_\phi^{\frac \eta {p}}xh_\phi^{\frac {1-\eta} {p}}\in L^{p}(\M)$ for $x\in\M,\ 0\le\eta\le1,$
\[
    \L_p^H(h_\phi^{\frac \eta {p}}xh_\phi^{\frac {1-\eta} {p}}):=h_\phi^{\frac \eta {p}}L(x)h_\phi^{\frac {1-\eta} {p}}.
\]
\end{definition}
\begin{definition}[Lindbladian Generator in Kosaki's $L^p$ Spaces]\label{def:lind-k-lp}
Under the assumptions of \Cref{def:lind-h-lp}, one can define the $L^p(\M,\phi)_\eta$ version
\[   \L_p^K: \Dom(\L_p^K) \subset L^p(\M,\phi)_\eta \to L^p(\M,\phi)_\eta\]

given by dense elements\footnote{Indeed, by \cite[Lemma~1.1]{Junge_Xu_2003}, $h_\phi^{\frac \eta p} \M h_\phi^{\frac{1-\eta}{p}}$ is dense in $L^p(\M)$ for $0<p<\infty, 0\le\eta\le1$.} $\displaystyle \iota_\eta(x):=h_\phi^\eta xh_\phi^{1-\eta}\in L^p(\M,\phi)_\eta$, for $x\in\M$
\[
    \L_p^K(\iota_\eta (x)):=\iota_\eta(L(x)).
\]
\end{definition}
We denote the semigroup generated by $\L_p^{\{H,K\}}$ as $(T_t)_p^{\{H,K\}} := e^{-t \L_p^{\{H,K\}}}$. The two definitions are related via the isometric isomorphism in \Cref{eq:h-k-iso}; specifically,
\[
    \iota_{\eta,p}\big( \L_p^H(h_\phi^{\frac{\eta}{p}} x h_\phi^{\frac{1-\eta}{p}}) \big)
    = \L_p^K\big( \iota_\eta(x) \big).
\]
Based on this identification, we omit the superscripts $\{H,K\}$ and write $\L_p$ when the context is clear.
We can also define the corresponding conditional expectations established in \cite[Section~2]{Junge_Xu_2003}.
\begin{definition}[Conditional expectation in $L^p$ spaces]\label{def:E_p}
    Under the assumptions of \Cref{def:Hlp}, let $\N\subset \M$ be a von Neumann subalgebra.
Suppose $E_\phi:\M\to \N$ is a $\phi$-preserving conditional expectation, i.e., $\phi\circ E_\phi=\phi$. Consider dense elements $h_\phi^{\frac \eta {p}}xh_\phi^{\frac {1-\eta} {p}}\in L^{p}(\M)$ for $x\in\M,\ 0\le\eta\le1,$ and the bimodule properties of $E_\phi$, the corresponding $L^p$ version $E_p^\phi$\footnote{By \cite[Lemma~2.1]{Junge_Xu_2003} for any $0\leq\eta_1,\eta_2\leq 1$, the $L^p$ conditional expectation are identical, thus $E_p^\phi$ have no $\eta$ dependence.} is:
    \begin{equation}\label{eq:Ep}
        E_p^\phi(h_\phi^{\frac \eta {p}}xh_\phi^{\frac {1-\eta} {p}}):=h_{\phi|_\N}^{\frac \eta {p}}E^\phi(x)h_{\phi|_\N}^{\frac {1-\eta} {p}}
    \end{equation}
\end{definition}
It is also known that $L^p(\N)=L^p(\N,\phi|_\N)$ can be naturally isometrically identified with a subspace of $L^p(\M)$. From now on we identify $\phi$, $\iota_\eta$, $\sigma_t$, and $\tau$
(the canonical trace from \cref{eq:canonical-trace}) with their restrictions on $\N$,
and do not distinguish them unless needed for clarity.
Under these identifications, we summarize the compatible QMSs in \Cref{fig:int}.
\begin{figure}[htbp]
  \centering
\begin{tikzcd}[row sep=3em,column sep=5em]
\mathcal{M}
  \arrow[d,"T_t"]
  \arrow[r,"\iota_\eta"]
&  L^1(\M)
  \arrow[d,"(T_t)_*"']
& L^p(\mathcal{M},\phi)_\eta \cong L^p(\M)
  \arrow[l,hook']
  \arrow[d,"(T_t)_p"]
\\
\mathcal{M}
  \arrow[r,"\iota_\eta"]
  \arrow[d,"E"]
&  L^1(\M)
\arrow[d,"E_1"]
& L^p(\mathcal{M},\phi)_\eta\cong L^p(\M)
  \arrow[l,hook']
  \arrow[d,"E_p^\phi"]
  \\
  \N
  \arrow[r,"\iota_\eta"]
  &
 L^1(\N)
  &
L^p(\N,\phi)_\eta\cong L^p(\N)    \arrow[l,hook']
\end{tikzcd}
\caption{Kosaki's $L^p$-interpolation spaces for $\mathcal{M}$ and their compatible QMSs, via the functoriality of complex interpolation \cite{Bergh_Löfström_1976}.}
  \label{fig:int}
\end{figure}
\begin{lemma}[\text{\cite[Lemma~2.2]{Junge_Xu_2003}}]\label{lem:norm-eq}
Under the assumptions of \Cref{def:E_p}, for each $1\le p<\infty$, $E_p^\phi$ extends to a contractive projection, hence isometric on its range.

In particular, in Haagerup's $L^p$ spaces
\[
\|E_p^\phi(x)\|_{L^p(\M)} \le \|x\|_{L^p(\N)}, \quad x\in L^p(\M),
\]
with equality for $x\in L^p(\N)$.
Equivalently, in Kosaki's $L^p$ spaces,
\[
\|\iota^{\phi|_\N}_{\eta,p}(E^\phi(x))\|_{p,\phi,\eta}=\|E_1^\phi(\iota^\phi_{\eta,p}(x))\|_{p,\phi,\eta}\le \|\iota^{\phi|_\N}_\eta (x)\|_{p,\phi|_\N,\eta}, \quad x\in \M,
\]
with equality for $x\in \N$.
\end{lemma}
We next explore the compatible definitions for gradient forms in both $L^p$ spaces and their identification.
\begin{definition}[Gradient form in Haagerup's $L^p $ spaces]\label{def:GammaLp-h}
Let
\[
\mathfrak D_p
:=
\left\{(x,y):
x,y\in\Dom(\mathcal L_p),\
x^*y\in\Dom(\mathcal L_{p/2})
\right\}.
\]
The Haagerup $L^p$ gradient form for the Lindbladian $\L_p$ is the sesquilinear map
\[
\Gamma_p:\mathfrak D_p\longrightarrow L^{p/2}(\M)
\]
defined by
\begin{equation}\label{eq:gf-h}
\Gamma_p(x,y)
=
\frac12\left(
\mathcal L_p(x^*)y
+x^*\mathcal L_p(y)
-\mathcal L_{p/2}(x^*y)
\right).
\end{equation}
For $(a,a)\in\mathfrak D_p$, we use the shorthand
\[
\|a\|_{\Gamma_p}
:=
\|\Gamma_p(a,a)\|_{L^{p/2}(\M)}^{1/2}.
\]
The same convention applies when the $L^p$ gradient form carries additional
decorations, such as $\Gamma_{L_\varepsilon,p}$ or $\Gamma_p^k$.
\end{definition}
\begin{definition}[Gradient form in Kosaki's $L^p$ spaces]\label{def:GammaLp-k}
For $x,y\in \M_{a}^{\sigma}\cap\Dom(L)$, define the gradient form compatible with the Kosaki interpolation space by
\[
    \Gamma^{(p)}_{\eta}(x,y):= \Gamma(\sigma_{- \frac {i\eta} {p}}(x),\sigma_{- \frac {i\eta} {p}}(y)),
\]
Here $\M_{a}^{\sigma}$ is the set of all entire elements of $\M$\footnote{An element $x\in\M$ is said to be entire if the function $\R\ni t\mapsto \sigma_t(x)\in\M$ can be extended to an $\M$-valued entire function over $\mathbb{C}$. It follows from \cite[Lemma~2.3]{takesaki_theory_2003-2} that $\M_a^\sigma$ is a $\sigma$-weakly dense $*$-subalgebra of $\M$.}, $\sigma_{- \frac {i\eta} {p}}(x)=h_\phi^{\frac{\eta}{p}}xh_\phi^{-\frac{\eta}{p}}$ is the modular automorphism at complex time $t=\frac{-i\eta}{p}$, and $\Gamma$ is defined in \cref{eq:Gamma}.
For such $x$, we use the shorthand
\[
\|x\|_{\Gamma^{(p)}_\eta}
:=
\bigl\|\iota_{\frac12}\bigl(\Gamma^{(p)}_\eta(x,x)\bigr)\bigr\|_{\frac p2,\phi,\frac12}^{1/2}.
\]
\end{definition}
The identification between two gradient forms is for dense elements,
\[
\Gamma_{p}\big(h_\phi^{\frac{\eta}{p}} x h_\phi^{\frac{1-\eta}{p}},  h_\phi^{\frac{\eta}{p}} y h_\phi^{\frac{1-\eta}{p}}\big)
= h_\phi^{\frac{1}{p}}   \Gamma_{\eta}^{(p)}(x, y)   h_\phi^{\frac{1}{p}}, \quad  x,y\in\M
\]
In the tracial case, one can replace $h_\phi$ by $D_\phi$, where $D_\phi$ is the density of $\phi$ with respect to the trace, i.e., $\phi(\cdot) = \tau(D_\phi \cdot)$.
\subsection{Translation towards non-trace-symmetric QMS}
In this section, we rewrite \Cref{main:tracial} and \Cref{cor:non-self-adjoint-ppp} for $\tau$-symmetric QMS in $(\M,\tau)$ in terms of relevant $L^p$ quantities as in \cref{sec:lind-lp} with an additional assumption of the KMS-$\phi$-detailed balance condition. The state $\phi$ admits a density $D_\phi$ with respect to the trace $\tau$, satisfying $\phi(\cdot) = \tau(D_\phi \cdot)$.
 The following results follow directly from rewriting and the factorization property in \cref{eq:kms-2}.

\begin{proposition}
     [Haagerup's $\Poin (p,p)$]\label{prop:h-lp-conj}
    Given the same conditions as the theorem \ref{main:tracial} and assuming further that the QMS $(T_t)_{t\ge0}$ is KMS-$\phi$-detailed balance, the following inequality holds in Haagerup's $L^p$ spaces for all $p = 2$ and $p \ge 3$, and for all dense $x\in\M_{sa}$, $a:=D_\phi^{\frac 1{2p}} xD_\phi^{\frac 1{2p}}  \in L^p(\M)_{{sa}}$ in the domain:
\[
\|a-E_p^\phi(a)\|_{L^p(\M)} \leq \frac{p}{\sqrt{2\alpha}}\|a\|_{\Gamma_p}.
\]
Furthermore, for $a\in L^p(\M)$ in the domain,
\[
\|a-E_p^\phi(a)\|_{L^p(\M)} \leq \frac{p}{\sqrt{2\alpha}}
\bigl(\|a\|_{\Gamma_p}+\|a^*\|_{\Gamma_p}\bigr).
\]
\end{proposition}
\begin{proposition}
    [Kosaki's $\Poin (p,p)$]\label{prop:k-lp-conj}
       Given the same conditions as in \Cref{prop:h-lp-conj}, for every $p=2$ and $p\geq 3$ and $x\in\M_{sa}$ in the domain, we have the inequality in Kosaki's $L^p$ spaces:
\[
  \bigl\|\iota_{\frac12}(x-E(x))\bigr\|_{p,\phi,\frac{1}{2}}
   \le
  \frac{p}{\sqrt{2\alpha}}
  \|x\|_{\Gamma^{(p)}_{\frac12}}.
\]
Furthermore, for $x\in\M$ in the domain,
\[
  \bigl\|\iota_{\eta}(x-E(x))\bigr\|_{p,\phi,\eta}
   \le
  \frac{p}{\sqrt{2\alpha}}
  \Bigl(
      \|x\|_{\Gamma^{(p)}_\eta}
      +\|x^*\|_{\Gamma^{(p)}_{1-\eta}}
  \Bigr).
\]
In particular, the one-sided version is
\begin{equation}\label{eq:non-self-adjoint}
  \bigl\|(x-E(x)) D_{\phi}\bigr\|_{p,\phi,0}
   \le
  \frac{p}{\sqrt{2\alpha}}
  \Bigl(
      \|x\|_{\Gamma^{(p)}_0}
      +\|x^*\|_{\Gamma^{(p)}_1}
  \Bigr),
\end{equation}
\end{proposition}
Moreover, \Cref{prop:h-lp-conj,prop:k-lp-conj} remain well-defined for non-tracial von Neumann algebras: one needs only replace the canonical density $D_{\phi}$ by $h_{\phi}$ (cf.\ item (v) of \Cref{prop:h-lp}). The replaced $h_\phi$ is defined in the absence of a trace.

As discussed in \Cref{rmk:gns-kms}, when we have a $\sigma$-finite $(\M,\phi)$, the $\tau$-symmetric condition is no longer available. Therefore, the natural replacement is the stronger GNS--$\phi$--detailed balance. We formulate our main theorem below and prove it in later sections.
\begin{theorem}[$\Poin(p,p)$ for GNS-$\phi$-detailed balance]\label{thm:GNS-Lp-ppp}
Let $(\M,\phi)$ be a $\sigma$-finite von Neumann algebra equipped with a faithful normal state $\phi$, and let $(T_t)_{t\ge0}$ be a GNS-$\phi$-detailed-balance QMS with fixed-point algebra $\N$. Suppose $(T_t)_{t\ge0}$ has a spectral gap \[\alpha:=  \inf_{\substack{x\in\ker E\\ x\neq 0}}
  \frac{\phi\bigl(x^*L(x)\bigr)}
       {\phi(x^*x)}
  >0.\]
Let $E_p^\phi$ be the induced conditional expectation from $L^p(\M)$ onto $L^p(\N)$. Then, for $p=2$ or $p\ge 3$ and $a\in L^p(\M)_{sa}$ in the domain, the following holds in the Haagerup $L^p$-space:
\[
\|a-E_p^\phi(a)\|_{L^p(\M)} \leq \frac{p}{\sqrt{2\alpha}}\|a\|_{\Gamma_p}.
\]
Equivalently, using symmetric Kosaki's $L^p$ spaces notation, for $x\in \M_{sa}$ in the domain
\[
  \bigl\|\iota_{\frac12}(x-E^\phi(x))\bigr\|_{p,\phi,\frac12}
   \le
  \frac{p}{\sqrt{2\alpha}}
  \|x\|_{\Gamma^{(p)}_{\frac12}}.
\]
\end{theorem}
\begin{cor}[GNS-$\phi$-detailed balance $\Poin(p,p)$ for non-self-adjoint elements]\label{cor:non-self-adjoint}
Under the construction of \Cref{thm:GNS-Lp-ppp}, for $a\in L^p(\M)$ in the domain,
\[
\|a-E_p^\phi(a)\|_{L^p(\M)} \leq \frac{p}{\sqrt{2\alpha}}
\bigl(\|a\|_{\Gamma_p}+\|a^*\|_{\Gamma_p}\bigr).
\]
Equivalently, using Kosaki's $L^p$ spaces notation, for $\eta\in[0,1]$, $x\in \M$ in the domain
\[
  \bigl\|\iota_{\eta}(x-E(x))\bigr\|_{p,\phi,\eta}
   \le
  \frac{p}{\sqrt{2\alpha}}
  \Bigl(
      \|x\|_{\Gamma^{(p)}_\eta}
      +\|x^*\|_{\Gamma^{(p)}_{1-\eta}}
  \Bigr).
\]
\end{cor}
\begin{proof}
As in \Cref{cor:non-self-adjoint-ppp}, write $a=\Re(a)+i \Im(a)$ and apply the bound to $\Re(a)$ and $\Im(a)$.
\end{proof}
\begin{remark}[Domain convention]
In the tracial setting, Cipriani and Sauvageot \cite{cipriani_derivations_2003} identified and extended the domain of the gradient form. Indeed, they showed that $\Dom(L^{1/2})$ is a *-algebra by constructing
a derivation such that
\[
\Dom(\overline\delta)=\Dom(L^{1/2}),
\qquad
\|\overline\delta(x)\|_{\mathcal H}^{2}
=
\|L^{1/2}x\|_2^{2}=\|\Gamma(x,x)\|_1.
\]
Moreover, if the state is tracial, then
$T_t(x)\in \Dom(L^{1/2})\cap\Dom(L)$ for every $x\in\M$ and $t>0$.
Thus, $x_n=T_{1/n}(x)$ gives an approximating sequence.
Although the gradient form is initially defined only on a regular core, one would like to define a gradient seminorm for general elements of $L^p(\M)$ by approximation.
For a general $a\in L^p(\M)$, the gradient seminorm in
\Cref{thm:GNS-Lp-ppp} is therefore understood via analytic approximation
\begin{equation}\label{eq:Gamma-p-norm-approximation}
\|a\|_{\Gamma_p}
:=
\inf_{\substack{
 a_n\in\mathcal C_p\\
 \|a_n-a\|_{L^p(\M)}\to0}}
\liminf_{n\to\infty}
\|a_n\|_{\Gamma_p},
\end{equation}
where
\[
\mathcal C_p
:=
\left\{
h_\phi^{1/(2p)}x h_\phi^{1/(2p)}:
x\in\M_a^\sigma\cap\Dom(L)
\right\},
\]
and the infimum is $+\infty$ if no such approximating sequence exists.
We defer a detailed treatment of the domain issue to \Cref{app:conv}.
\end{remark}
\subsection{Haagerup's reduction}\label{sec:reduction}
To generalize our results to the non-tracial case, we introduce Haagerup's reduction \cite{haagerup2009reductionmethodnoncommutativelpspaces}, which shows that every von Neumann algebra admits finite approximations. We later show that the tracial $\Poin(p,p)$ inequality can be transferred to a general von Neumann algebra via Haagerup reduction.
Throughout this section, let $  G = \bigcup_{n\ge1}2^{-n}\mathbb{Z}  \subset \R$
and let $\M$ be a $\sigma$-finite\footnote{Mutually orthogonal projections in $\M$ are at most countable.}  von Neumann algebra equipped with a normal faithful state $\phi$.  We denote by
$\sigma_t^{\phi}$ the modular automorphism group of $\phi$ on $\M$, and form the crossed-product
$\wh{\M}:=\M\rtimes_{\sigma^\phi}G$. There is a canonical normal conditional expectation $\cE:\wh\M\to\M$\footnote{We avoid using the crossed product $\M\rtimes\R$ because, by \cite[Sec.~IX.4]{Takesaki_1979}, the existence of a normal conditional expectation onto $\M$ forces the factor type to be preserved.  In particular, even if $\M$ is type III, its crossed product with $\R$ is always semifinite (type II$_\infty$) due to the presence of a faithful trace.}.  Finally, let $\wh\phi=\phi\circ\cE$ be the dual weight on $\wh\M$, which is also a normal faithful state.

\begin{theorem}[Haagerup reduction {\cite[Theorem~2.1]{haagerup2009reductionmethodnoncommutativelpspaces}}]\label{thm:haagerup}
  Following the notation above, there exists an increasing sequence $(\M_k)_{k\ge1}$ of von Neumann subalgebras of $\wh\M$ satisfying:
  \begin{enumerate}[label=(\roman*)]
    \item each $\M_k$ is finite with normal faithful trace $\tau_k$;
    \item $\bigcup_{k\ge1}\M_k$ is w*-dense in $\wh{\M}$;
    \item there is a family of normal faithful conditional expectation $(\cE_k)_{k\ge1} \colon \wh\M  \to \M_k$
    such that
    \[
      \wh\phi\circ\cE_k = \wh\phi,
      \quad
      \sigma_t^{\wh\phi}\circ\cE_k = \cE_k\circ\sigma_t^{\wh\phi},
      \quad
      t\in\R;
    \]
    \item let $\phi_k(\cdot):=\wh{\phi}|_{\M_k}(\cdot)=\tau_k (d_k\cdot)$ for $d_k\in L(G)\subset \cZ(\wh{\M}_{\wh{\phi}})$ so that $\sigma_t^{\wh{\phi}}(d_k)=d_k$; Furthermore $a_k\le d_k\le a_k^{-1}$ for some $a_k\in\R^+$;
    \item $\lim_{k\to\infty}\cE_k(x)=x$ in $\sigma$-strong topology for any $x\in\wh\M$.
    \item (\cite[Section~4.2]{Gao_Junge_LaRacuente_Li_2025}) $\varphi_k \to \wh{\varphi}$ in the weak topology on $S(\M)$, for $\varphi\in S(\M)$ and $\varphi_k:=\wh{\varphi}|_{\M_k}$.
  \end{enumerate}
\end{theorem}

We apply Haagerup's reduction \Cref{thm:haagerup} to the GNS-$\phi$-detailed-balance QMS $(T_t)_{t\ge0}$. We call ${(\wh T_t)}_{t\ge0}$ its canonical extension to $\wh{\M}$, defined by $\wh T_t:=T_t\otimes\Id_{B(l_2(G))}|_{\wh \M}$. It is also shown in \cite[Section~4.2]{Gao_Junge_LaRacuente_Li_2025} to be a \emph{GNS-$\wh\phi$-detailed-balance} QMS with fixed-point algebra $\wh{\N}:=\N \rtimes_{\sigma^\phi}G$.

Moreover, the $\wh\phi$-preserving conditional expectation $\wh E:\wh{\M}\to\wh{\N}$
 is just the canonical extension of $E:\M\to\N$. It turns out that $(\wh{T_t})_{t\ge0}$ is compatible with the finite approximations in Haagerup reduction, summarized in the theorem below.
\begin{theorem}[\text{\cite[Lemma~4.5]{Gao_Junge_LaRacuente_Li_2025}}]\label{thm:fn-multiplicative-domain}
  Let $(T_t)_{t\ge0}\colon \M\to\M$ be a GNS-$\phi$-symmetric quantum Markov semigroup. With the notation above and $t\ge0$, the following hold:
  \begin{enumerate}[label=(\roman*)]
    \item For each $k$ and $s\in \mathbb R$,
\[
\widehat T_t\circ \sigma_s^{\phi_k}
=
\sigma_s^{\phi_k}\circ \widehat T_t .
\]
Consequently, $\widehat T_t$ commutes with $\mathcal E$, $\widehat{\mathcal E}$, and
$\mathcal E_k$. In particular, $\widehat T_t(\mathcal M_k)\subseteq \mathcal M_k$.
\item The restriction
      \[
        T_{t,k} := \wh T_t\bigm|_{\M_k}
         : \M_k \longrightarrow \M_k
      \]
      is a $\tau_k$-symmetric normal unital completely positive map.
    \item Let $\N_k:=\M_k\cap\wh\N\subseteq \M_k$ be the fixed-point algebra of $T_{t,k}$.  Then
      \[
        E_k := \wh E\bigm|_{\M_k}
         : \M_k \longrightarrow \N_k
      \]
      is the $\tau_k$-preserving conditional expectation onto $\N_k$.
  \end{enumerate}
\end{theorem}

Combining the preceding results with \cite[Eq.~(6.4)]{haagerup2009reductionmethodnoncommutativelpspaces}, which gives
\[
\wh E\circ\cE=\cE\circ\wh E
\quad\text{and}\quad
\wh E\circ\cE_k=\cE_k\circ\wh E,
\]
we obtain the following commuting diagrams \Cref{fig:haagerup-arrows}.
\begin{figure}[htbp]
  \centering
\begin{tikzcd}[row sep=2.5em,column sep=5em]
  \M
    \arrow[d,"T_t"']
  &
  \wh\M
    \arrow[d,"\wh T_t"]
    \arrow[l,"\cE"']
    \arrow[r,"\cE_k"]
  &
  \M_k
    \arrow[d,"T_{t,k}"]
\\
  \M
    \arrow[d,"E"']
  &
  \wh\M
    \arrow[d,"\wh E"]
    \arrow[l,"\cE"']
    \arrow[r,"\cE_k"]
  &
  \M_k
    \arrow[d,"E_k"]
\\
  \N
  &
  \wh\N
    \arrow[l,"\cE"']
    \arrow[r,"\cE_k"]
  &
  \N_k
\end{tikzcd}
\caption{Haagerup reduction commuting diagram: the quantum Markov semigroup $T_t$ on $\M$, its finite approximants $\M_k$, and the conditional expectations $E$, $E_k$ down to $\N$.}  \label{fig:haagerup-arrows}
\end{figure}

We have the following corollary on the tracial slices $(\M_k,\tau_k)$:
\begin{cor}\label{cor:d_k}
Using the notation above, for every $\alpha\in \R$ and every $x\in \M_k$,
\[
  T_{t,k}(d_k^\alpha x)
  =
  d_k^\alpha T_{t,k}(x),
\]
and similarly
\[
  T_{t,k}(x d_k^\alpha)
  =
  T_{t,k}(x)d_k^\alpha .
\]
\end{cor}

\begin{proof}
Since $d_k\in \M_k$ and $T_{t,k}=\widehat T_t|_{\M_k}$,
\[
  T_{t,k}(d_k)=\widehat T_t(d_k)=d_k .
\]
Because $T_{t,k}$ is unital completely positive and $\tau_k$-preserving,
Kadison--Schwarz and faithfulness of $\tau_k$ imply
\[
  T_{t,k}(d_k^2)=T_{t,k}(d_k)^2=d_k^2 .
\]
Thus $d_k$ lies in the multiplicative domain of $T_{t,k}$. By functional
calculus, the same holds for $d_k^\alpha$, and $T_{t,k}(d_k^\alpha)=d_k^\alpha$.
Hence
\[
  T_{t,k}(d_k^\alpha x)
  =
  T_{t,k}(d_k^\alpha)T_{t,k}(x)
  =
  d_k^\alpha T_{t,k}(x).
\]
The right multiplication identity is identical.
\end{proof}
\subsection{Proof of \texorpdfstring{$L^p$}{Lp} version of \texorpdfstring{$\Poin(p,p)$}{Poin(p,p)}}\label{sec:pf-lp-p}

We are now ready to prove our main theorem \Cref{thm:GNS-Lp-ppp}. Unless otherwise stated, we use the notation established in previous sections. The outline is to prove that $\PI(p,p)$ holds on the tracial slices $(\M_k,\tau_k)$ for all $k\ge1$ and extend it to $(\M,\phi)$ by Haagerup's reduction.
\begin{proof}
We first verify that the spectral gap $\alpha$ is inherited by $T_{t,k}$.
Combining this with (ii) of \Cref{thm:fn-multiplicative-domain},
it follows from \Cref{main:tracial} that $\Poin(p,p)$ holds on each slice $(\M_k,\tau_k)$.
Indeed, given $\alpha>0$ spectral gap for $T_t$, we have for $x\in\M$,
\[
    \|T_t(x)-E(x)\|_{L^2(\M)}\le e^{-\alpha t}\|x-E(x)\|_{L^2(\M)}.
\]
Recall that $\wh T_t:=T_t\otimes \Id_{B(\ell_2(G))}|_{\wh \M}$.
By direct calculation, or equivalently by \Cref{prop:Gap-tensor}, using that the spectral gap of $\Id$ is $\alpha_{\Id}=\infty$,
we see that $\wh T_t$ has the same spectral gap $\alpha$ as $T_t$. By isometry of norms for $\cE_k$ on $L^2$ spaces using \Cref{lem:norm-eq}, it follows that $\wh T_t$ has the same spectral gap $\alpha$ with respect to $\phi_k(\cdot)=\tau_k(d_k\cdot)$ when restricted to each slice $(\M_k,\tau_k)$ for $k\ge 1$.
Define
\[
U_k:L^2(\M_k,\phi_k)\to L^2(\M_k,\tau_k),
\qquad U_k(x)=xd_k^{1/2}.
\]
By (iv) of \Cref{thm:haagerup}, $d_k$ is boundedly
invertible, so $U_k$ is unitary. Moreover, since
$d_k^{1/2}\in\N_k$, \Cref{cor:d_k} and the
$\N_k$-bimodularity of $E_k$ imply
\[
U_kT_{t,k}=T_{t,k}U_k,
\qquad
U_kE_k=E_kU_k.
\]
Therefore $(T_{t,k},E_k)$ on $L^2(\M_k,\phi_k)$ is
unitarily equivalent to $(T_{t,k},E_k)$ on
$L^2(\M_k,\tau_k)$. In particular, $T_{t,k}$ has the
same spectral gap $\alpha$ with respect to $\tau_k$.

In the sequel, we adopt the convention that all objects associated with the restricted semigroup $T_{t,k}$ carry the subscript $k$.  In particular, we set
\[
L_k  := \wh L\bigl|_{\M_k},
\qquad
\L_{p,k}^{K}:=\wh{\L}^{K}_p\bigl|_{\M_k},
\qquad
\Gamma^{(p,k)}_{\eta}
 := \wh\Gamma^{(p)}_{\eta}\bigl|_{\M_k\times\M_k},
\]
so that $L_k$, $\L^K_{p,k}$ and $\Gamma^{(p,k)}_{\eta}$ are the Lindbladian, the Kosaki's $L^p$ Lindbladian and $L^p$ gradient form of $T_{t,k}$ on $\M_k$.
Given $k\ge1$, let $x_k\in (\M_{k})_{sa}\cap\wh\M_a^{\sigma}\cap \Dom(L_k)$, where $\wh \M_{a}^{\sigma}$ the set of all entire elements of $\wh \M$, and $d_k$ the density of $\phi_k:=\wh{\phi}|_{\M_k}$ under trace $\tau_k$. By (iv) of \Cref{thm:haagerup}, $d_k\in L(G)\subset\wh{\N}:=\N \rtimes_{\sigma^\phi}G$ and $d_k\in \M_k$, it follows that $d_k\in \N_k=\M_k\cap\wh{\N}$. We can assert the following
    \[
        \L_{p,k}^K(d_k^{\frac{1}{2}} x_kd_k^{\frac{1}{2}})= d_k^{\frac{1}{2}}L_k(x_k)d_k^{\frac{1}{2}},
    \]
which agrees with \Cref{def:lind-h-lp}. Thus, we can write the $\Poin(p,p)$ in Kosaki's $L^p$ spaces as in \Cref{prop:k-lp-conj}.
For all $p = 2$ and $p \ge 3$, and let $\iota_\eta(x)=h_\phi^{\eta} x h_\phi^{(1-\eta)}$ \footnote{We omit the $\phi$-subscript on $\iota$ whenever the norm's subscript makes the state clear.
It is also known that $L^p(\N)$ is naturally isometrically embedded in $L^p(\M)$
\cite{Junge_Xu_2003} whenever there exists a conditional expectation $E:\M\to\N$,
so no further distinction is needed.} for $\eta\in[0,1]$,

\[
  \bigl\|\iota_{\frac12}(x_k-E_k(x_k))\bigr\|_{p,\phi_k,\frac{1}{2}}
   \le
  \frac{p}{\sqrt{2\alpha}}
  \bigl\|\iota_{\frac12}  \bigl(\Gamma_{\frac 12}^{(p,k)}(x_k,x_k)\bigr)\bigr\|_{\frac{p}{2},\phi_k,\frac 12}^{1/2}.
\]
As we are working with self-adjoint elements we denote $\|\cdot \|_{p,\phi,\frac{1}{2}}$ just as $\|\cdot \|_{p,\phi}$ in the following.
Taking the limit $k\to\infty$, the above relations remain valid, and hence $\Poin(p,p)$ holds in $\wh{\M}$ and thus in $\M$.

To see this, for ${x}\in\cB:={\M}_{sa}\cap{\M}^\sigma_a\cap \Dom( L)\subset \M \subset \wh \M$, we have $x_k:=\cE_k({x})\in(\M_k)_{sa}$, $\cE_{k,1}\circ\iota_\eta=\iota_\eta\circ\cE_k$ and $E_k\circ\cE_k=\cE_k\circ \wh E=\cE_k\circ E$\footnote{Since $\wh E$ is the canonical extension of $E$, we write $E$ for $\wh E$ when acting on $x\in\M$.} from
\Cref{fig:int,fig:haagerup-arrows}, the left-hand side rewrites as
\begin{equation}\label{eq:1}
     \bigl\|\iota_{\frac12}(x_k-E_k(x_k))\bigr\|_{p,\phi_k}=\bigl\|\iota_{\frac12}(\cE_k(x)-E_k(\cE_k(x)))\bigr\|_{p,\phi_k}=\bigl\|\cE_{k,1}\bigl(\iota_{\frac12}(x- E(x))\bigr)\bigr\|_{p,\phi_k}.
\end{equation}
Taking the limit $k\to \infty$ and applying (v) of \Cref{thm:haagerup}, $\lim_{k\to\infty}\cE_k(x)=x$ in $\sigma$-strong topology for any $x\in\cB\subset\wh\M$ give us
\begin{equation}\label{eq:2} \lim_{k\to\infty}\bigl\|\cE_{k,1}\bigl(\iota_{\frac12}(x- E(x))\bigr)\bigr\|_{p,\phi_k}=\bigl\|\iota_{\frac12}(x- E(x))\bigr\|_{p,\wh \phi}=\bigl\|\iota_{\frac12}(x- E(x))\bigr\|_{p,\phi},
\end{equation}
where the last equality uses \Cref{lem:norm-eq} where the $L^p$ extensions of $\cE:\wh\M\to\M$ is a contractive projection.

Similarly, for the RHS, we can remove one conditional expectation as in \cref{eq:1,eq:2},
\[ \|\iota_{\frac{1}{2}}\bigl(\Gamma_{\frac 12}^{(p,k)}(x_k,x_k)\bigr)\bigr\|_{\frac{p}{2},\phi_k }^{1/2}\le\|\iota_{\frac{1}{2}}\bigl(\wh\Gamma_{\frac 12}^{(p)}(x,\cE_k(x))\bigr)\bigr\|_{\frac{p}{2},\wh \phi }^{1/2}\]

We know $\wh\Gamma_{\frac 12}^{(p)}(x,\cE_k(x))$ is linear on $\cE_k(x)$ from \Cref{def:GammaLp-k},
\[
    \wh\Gamma_{\frac 12}^{(p)}(x,\cE_k(x))=\Gamma(\sigma^{\wh\phi}_{- \frac {i} {2p}}(x),\sigma^{\wh\phi}_{- \frac {i} {2p}}(\cE_k(x)))=\Gamma(\sigma^{\wh\phi}_{- \frac {i} {2p}}(x),\cE_k(\sigma^{\wh\phi}_{- \frac {i} {2p}}(x))),
\]
where the second equality is due to $      \sigma_t^{\wh\phi}\circ\cE_k = \cE_k\circ\sigma_t^{\wh\phi}$ by (iii) of \Cref{thm:haagerup}. Let $y:=\sigma^{\wh\phi}_{- \frac {i} {2p}}(x)=\sigma^{\wh\phi}_{- \frac {i} {2p}}\bigl|_{\M}(x)=\sigma^{\phi}_{- \frac {i} {2p}}(x)\in\M$ for $x\in \M_a^\sigma$, and $y_k:=\cE_k(y)$. We aim to prove that
\begin{equation}\label{eq:gamma-pnorm}
    \lim_{k\to\infty}\|\iota_{\frac{1}{2}}\bigl(\wh\Gamma(y,y_k)\bigr)\bigr\|_{\frac{p}{2},\wh \phi,\frac 12}^{1/2}\le \|\iota_{\frac{1}{2}}\bigl(\wh\Gamma(y,y)\bigr)\bigr\|_{\frac{p}{2},\wh \phi,\frac 12}^{1/2}=\|\iota_{\frac{1}{2}}\bigl(\Gamma(y,y)\bigr)\bigr\|_{\frac{p}{2}, \phi,\frac 12}^{1/2},
\end{equation}
where the last equality is due to the isometry of the norm from \Cref{lem:norm-eq}.
There are technical complications due to the domain issue for the gradient form $\Gamma$: $\Dom(L)$ is not closed under multiplication.
In particular, for $y,y_k\in\Dom( L)$ one may have
$y^*y_k\notin\Dom(L)\subset \Dom({\wh L})$. To work around this issue, we employ the
$\epsilon$-regularization, so that all
expressions are well defined and the domain problem does not obstruct the passage to the limit. Recall in this framework, for all $x,y\in\cA:=\Dom( L)\cap \M\subset\Dom( \wh L)\cap \wh \M$,
\begin{equation}\label{eq:regulated-Gamma}
        \wh \Gamma_\epsilon(x,y)
 :=
\frac{1}{2}\Bigl(
x^*\wh L_\epsilon(y)+\wh L_\epsilon(x)^*y-\wh L_\epsilon(x^*y)
\Bigr), \text{ where } \wh L_\epsilon:=\frac{\wh L}{1+\epsilon \wh L}.
\end{equation}

It follows from \cite[Prop.~2.5]{cipriani_derivations_2003} that, for $\epsilon>0$, the resolvent maps $(I+\epsilon \wh{L})^{-1}$ are bounded, completely positive contractions on $L^2({\wh\M})$ and the regulated operators $ \wh L_\epsilon$
are completely bounded endomorphisms of $\wh \M$.

These properties guarantee that the regularized gradient form above is
well defined and uniformly bounded (\Cref{rmk:unif-bdd}).
As shown in \Cref{app:regular}, detailed balance for $L$ passes to the
regulator $L_\epsilon$, so that the tracial
$\PI(p,p)$ inequality holds for the regularization $ L_{k,\epsilon}={L_k}({I+\epsilon  L_k})^{-1}$ on $\M_k$, and the spectral gap satisfies
$\alpha_\epsilon ={\alpha}/ ({1+\epsilon\alpha})$.

 Meanwhile, we restrict to the simplified case in which we assume that $\A$ is a weak$^*$-dense $*$-algebra, given by $\A = \Dom(L)\cap \M$. By \Cref{app:dom-strong}, $\wh\Gamma_\epsilon(x,y)\to \wh\Gamma(x,y)$ strongly for all $x,y\in\A$.
To prove \cref{eq:gamma-pnorm}, we first work with the regularized operators $\wh\Gamma_\epsilon$. By (v) of \Cref{thm:haagerup}, we have $\lim_{k\to\infty} y_k = y$ in the $\sigma$-strong topology. Since
$
\wh L_\epsilon = \wh L(1+\epsilon \wh L)^{-1}
$
is bounded for each $\epsilon>0$, it is straightforward to check that
\begin{equation}\label{eq:strong-gamma}
      \lim_{k\to\infty}\wh\Gamma_\epsilon(y,y_k)
    = \lim_{k\to\infty}\frac12\bigl(
        \wh L_\epsilon(y^*)  y_k
        + y^* \wh L_\epsilon(y_k)
        - \wh L_\epsilon(y^* y_k)
    \bigr)
    = \wh\Gamma_\epsilon(y,y)
\end{equation}
strongly.
After taking $\epsilon\downarrow 0$ and applying \Cref{lem:p-norm-convergence}, we obtain the equality form of \cref{eq:gamma-pnorm} in $\frac{p}{2}$-norm. For the general case when there is no $*$-algebra $\A$, we use the weak$^*$ convergence definition extended by \cite[Lemma~9.1]{cipriani_derivations_2003} (\cref{eq:Gamma-epsilon}).
\[
\Gamma(x,y):=\mathrm{w}^*\text{-}   \lim_{\epsilon\downarrow0}\Gamma_\epsilon(x,y).
\]

For an arbitrary self-adjoint $x\in\M_{sa}$, we consider a family of analytic approximations $\psi_m(x)\in\M_a^\sigma$.
Recall from \cite[Chap.~VIII.2.3]{takesaki_theory_2003-2} that these approximations may be defined by
\[
\psi_m(x):=\int_\R \sigma_t^\phi(x)K_m(t)\,dt,
\]
where $K_m$ is a family of suitable Gaussian kernels.
Each map $\psi_m$ commutes with the GNS-detailed-balance semigroup
$\{T_t\}_{t>0}$ because the semigroup commutes with the modular automorphism group $\{\sigma_t^\phi\}_{t\in\R}$.
Consequently, $T_{1/n}(\psi_m(x))\in\cB$.
Moreover, $T_{1/n}(\psi_m(x))\to x$ strongly as $n,m\to\infty$.
We apply the $\Poin{(p,p)}$ inequality to these approximating elements and pass to the limit on both sides, interpreting the $\Gamma_p$-norm from \Cref{eq:Gamma-p-norm-approximation}.

By the density of $\bigcup_k\M_k$, stated in part~(ii) of \Cref{thm:haagerup}, the above arguments establish the $L^p$ version of $\Poin(p,p)$ on $(\M,\phi)$ for every $x\in\M_{sa}\cap\Dom(L)$.
Equivalently, in Haagerup's functorial formulation of $L^p$-spaces, for every
$a=h_\phi^{1/(2p)}x h_\phi^{1/(2p)}\in L^p(\M)_{sa}$ with
$x\in\M_{sa}\cap\Dom(L)$, we have
\[
\|a-E_p^\phi(a)\|_{L^p(\M)}
\leq
\frac{p}{\sqrt{2\alpha}}
\|\Gamma_p(a,a)\|_{L^{p/2}(\M)}^{1/2}.
\]
Since Haagerup's functorial construction is independent of any particular choice of trace, this $L^p$ inequality is canonically well-defined for all $a\in L^p(\M)_{sa}$.
We defer the details of the general case to \Cref{app:conv}.
\end{proof}

\section{Examples and Applications}\label{sec:examples}

\subsection{Khintchine inequalities with Rademacher random variables}
We present a degree-one example in which the $\Poin(p,p)$ inequality yields
Khintchine-type estimates. This example also shows that, in the non-self-adjoint
extension $\Poin(p,p)$ from \Cref{cor:non-self-adjoint}, the second term
(the $\|a^*\|_{\Gamma_p}$ contribution) is genuinely necessary and cannot be removed.

Let $\Omega=\{-1,1\}^n$ be equipped with the uniform probability measure, and let
$\varepsilon_1,\dots,\varepsilon_n$ denote the coordinate (Rademacher) functions.
We consider the von Neumann algebra
\[
\A := L_\infty(\Omega) \bar\otimes \M  \cong  L_\infty(\Omega,\M),
\]
equipped with the product state $\psi:=\mathbb E\otimes\phi$, where $\mathbb E$ is the uniform
expectation on $L_\infty(\Omega)$. Concretely, elements of $\A$ are bounded
$\M$-valued functions on $\Omega$; in particular, we will often restrict to the
degree-one subspace consisting of elements of the form
\[
x=\sum_{i=1}^n \varepsilon_i \times x_i,\qquad x_i\in\M.
\]

For each $i$, let $\bE_i:\A\to\A$ denote the conditional expectation onto the
von Neumann subalgebra generated by $\{\varepsilon_j:j\neq i\}$, acting trivially
on $\M$. We define the Lindbladian by
\[
    L=\sum_{i=1}^n L_i,\qquad L_i=\Id-\bE_i.
\]
One checks that the spectral gap of $L$ is $1$ and $L$ is GNS-$\psi$-detailed balance, since each $L_i$ is a projection; moreover $E_{fix}=\bE\otimes\Id_\M$ and $E_{fix}(x)=0$.
The associated gradient form for $x=\sum_{i=1}^n \varepsilon_i\times x_i$  decomposes as
\[
\Gamma=\sum_{i=1}^n \Gamma_i,\quad
\Gamma_i(x,x)= x_i^*x_i,
\quad
\Gamma_i(x^*,x^*)= x_i x_i^*.
\]
Denote by $h_\phi$ the density of $\phi$, and set
\[
a=xh_\phi^{1/p}=\sum_{i=1}^n \varepsilon_i\times a_i,
\qquad a_i=x_i h_\phi^{1/p}.
\]
Following \Cref{eq:Ep}, $E_p^\phi(\cdot)=E(\cdot)h_{\phi}^{\frac 1p}$ and thus $E_p^\phi(a)=0$. By \Cref{cor:non-self-adjoint} and density, for all $a_i\in L^p(\M)$ and $p=2$ or $p\ge3$,
\[
\begin{aligned}
\|a-E_p^\phi(a)\|_{L^p(\A)}
&=\Bigl\|\sum_{i=1}^n \varepsilon_i\times a_i\Bigr\|_{L^p(\A)} \\
&\le
\frac{p}{\sqrt{2}}
\Bigl\|(\sum_{i=1}^n a_i^*a_i)^{1/2}\Bigr\|_{L^p(\M)} \\
&\quad+
\frac{p}{\sqrt{2}}
\Bigl\|(\sum_{i=1}^n a_i a_i^*)^{1/2}\Bigr\|_{L^p(\M)}.
\end{aligned}
\]
Our application gives a suboptimal constant for noncommutative Khintchine inequalities in \cite{Lust-Piquard_Pisier_1991} with constant $O(\sqrt{p})$.
\subsection{Sub-exponential concentration inequalities}
It is well known that Poincar\'e inequalities imply concentration inequalities. We establish that our $\Poin(p,p)$ implies a sub-exponential concentration inequality, which in the semi-commutative case is established in \cite[Thm~2.7]{huang2021poincare}.

To obtain the desired concentration inequality, we need to resolve the
$p$-dependence of $\Gamma_p$, as a direct estimate is not immediate. We therefore return to the derivation framework of
\cite{JungeRicardShlyakhtenko2025}. The following lemma is stated for a
general tracial von Neumann algebra and will be applied below to each
tracial slice $(\M_k,\tau_k)$.

We summarize the noncommutative derivation framework of
\cite[Prop.~2.6, Prop.~2.22]{JungeRicardShlyakhtenko2025} as below.

\begin{lemma}[Noncommutative derivation]\label{lem:diff-quotient-derivation}
Let $(\M,\tau)$ be a tracial von Neumann algebra and let
\[
L_\varepsilon=\frac{I-T_\varepsilon}{\varepsilon}
\]
be the $\varepsilon$-regularized generator associated with a trace-symmetric
QMS. There exists a free extension
\[
\M\subset \M\otimes_{T_\varepsilon}\M \subset \mathscr M
\]
associated with the Hilbert module construction \cite{lance1995hilbert} of $T_\varepsilon$, and a
self-adjoint element $\xi_\varepsilon\in \mathscr M$, such that
\[
\partial_\varepsilon(x):=[x,\xi_\varepsilon],
\qquad x\in\M,
\]
defines a derivation into $L^2(\mathscr M)$.
Then the following holds.

\begin{enumerate}
\item[(i)] For $x,y\in \M$,
\[
2\Gamma_{L_\varepsilon}(x,y)
=
E_{\M}\bigl(
(\partial_\varepsilon(x))^*
\partial_\varepsilon(y)
\bigr).
\]
\item[(ii)] For $2\le p\le \infty$, there exist universal constants
$C>0$ such that for self-adjoint $x\in \M_{sa}$
\[
\sqrt{2}
\|\Gamma_{L_\varepsilon}(x,x)\|_{L^{p/2}(\M)}^{1/2}
\le
\|\partial_\varepsilon(x)\|_{L^p(\mathscr M)}
\le
C
\|\Gamma_{L_\varepsilon}(x,x)\|_{L^{p/2}(\M)}^{1/2}.
\]
In particular,
\begin{equation}\label{eq:infty-estimate}    \|\partial_\varepsilon(x)\|_{\infty}
\le
2\sqrt{2}\|\Gamma_{L_\varepsilon}(x,x)\|_\infty^{1/2}.
\end{equation}
\end{enumerate}
\end{lemma}
On the tracial slices $(\M_k,\tau_k)$, \Cref{cor:d_k} gives for $\alpha\ge0$, $
  T_{t,k}(d_k^\alpha x)
  =
  d_k^\alpha T_{t,k}(x),$ where $d_k\in \M_k$ is the density of $\phi_k:=\wh\phi|_{\M_k}$ with respect to $\tau_k$. As $T_{t,k}$ is
$\tau_k$-symmetric, the $L^p$-version identifies with $T_{t,k}$ on
self-adjoint elements
\[
  a=d_k^{1/2p}xd_k^{1/2p}\in L^p(\M_k),
  \qquad x\in(\M_k)_{\mathrm{sa}}.
\]
The following corollary is also for a tracial von Neumann algebra
$(\M,\tau)$, where $L^p(\M)$ is identified with $\M$. Essentially,
we apply the framework for $(\M_k,\tau_k)$ and drop the subscript $k$
to be notationally clear.

\begin{cor}\label{cor:p-derivation}
Let $(\M,\tau)$ be finite and tracial, and let $d\in \M$ be the density of
a normal faithful state $\phi$ with respect to $\tau$. Assume that $T_\varepsilon$
is $\tau$-symmetric and
\[
  T_\varepsilon(d^\alpha x)=d^\alpha T_\varepsilon(x),
  \qquad \alpha\ge0,\ \varepsilon>0,\ x\in\M .
\]
Following the notation of \Cref{lem:diff-quotient-derivation}, one has
\[
  [\xi_\varepsilon,d^\alpha]=0,
  \qquad \alpha\ge0.
\]
Consequently, the derivation associated with the $L^p$-version of the
Lindbladian $L_{\varepsilon,p}$ from \Cref{def:lind-h-lp} on dense elements is given by
\begin{equation}\label{eq:p-derivation}
      \partial_{\varepsilon}(a)
  :=
  [a,\xi_\varepsilon]
  =
  d^{1/2p}\partial_\varepsilon(x)d^{1/2p},
  \qquad
  a=d^{1/2p}x d^{1/2p},\ x\in \M_{\mathrm{sa}}.
\end{equation}
The $L^p$-gradient form \cref{eq:gf-h}, by identification of
$L_{\varepsilon,p}$ with $L_\varepsilon$ on dense elements and item~(i) of
\Cref{lem:diff-quotient-derivation}, satisfies
\[
  2\Gamma_{L_{\varepsilon,p}}(a,a)
  =
  E_{\M}\!\left(
    \partial_{\varepsilon}(a)^*
    \partial_{\varepsilon}(a)
  \right).
\]
\end{cor}

\begin{proof}
We use the free construction of $\mathscr M$ from
\cite[Section~2]{JungeRicardShlyakhtenko2025}. In this construction
\[
  \xi_\varepsilon=s(1_{[0,1/\varepsilon]}),
  \qquad s(f)=\ell(f)+\ell(f)^*,
\]
where $\ell(f)$ is the left creation operator used in the free
construction of $\mathscr M$.

By the multiplicative-domain property of $d^\alpha$,
\[
  T_\varepsilon(d^\alpha x)=d^\alpha T_\varepsilon(x),
  \qquad
  T_\varepsilon(xd^\alpha)=T_\varepsilon(x)d^\alpha .
\]
Thus, in $\M\otimes_{T_\varepsilon}\M$, the left and right actions of
$d^\alpha$ agree on the distinguished vector $\zeta=1\otimes1$:
\[
  d^\alpha\zeta=\zeta d^\alpha .
\]
Therefore, for every $f\in L^2(0,\infty)$,
\[
  d^\alpha(f\otimes\zeta)=(f\otimes\zeta)d^\alpha .
\]
Since $\ell(f)$ is the left creation operator
\[
\ell(f)(\zeta_1\otimes\cdots\otimes\zeta_n)
=
(f\otimes\zeta)\otimes \zeta_1\otimes\cdots\otimes\zeta_n,
\]
the relation $d^\alpha(f\otimes\zeta)=(f\otimes\zeta)d^\alpha$ implies
\[
  [\ell(f),d^\alpha]=0.
\]
Taking the adjoint gives the same commutation relation for $\ell(f)^*$.
Choosing $f=1_{[0,1/\varepsilon]}$, we get
\[
  [\xi_\varepsilon,d^\alpha]=0.
\]

Now let $a=d^{1/2p}xd^{1/2p}$. Since $\xi_\varepsilon$ commutes with
$d^{1/2p}$,
\[
  [a,\xi_\varepsilon]
  =
  d^{1/2p}[x,\xi_\varepsilon]d^{1/2p}
  =
  d^{1/2p}\partial_\varepsilon(x)d^{1/2p}.
\]
This gives the formula for $\partial_{\varepsilon}(a)$. The identity
for $\Gamma_{L_{\varepsilon,p}}$ follows from item~(i) of
\Cref{lem:diff-quotient-derivation} applied to the $L^p$-gradient form
\cref{eq:gf-h}:
\[
  2\Gamma_{L_{\varepsilon,p}}(a,a)
  =
  E_{\M}\!\left(
    [a,\xi_\varepsilon]^*[a,\xi_\varepsilon]
  \right)
  =
  E_{\M}\!\left(
    \partial_{\varepsilon}(a)^*
    \partial_{\varepsilon}(a)
  \right).
\]
\end{proof}
We can establish the following Corollary for $\Poin(p,\infty)$.
\begin{cor}[$\Poin(p,\infty)$ with $O(p)$ constant]\label[corollary]{cor:lip-norm-estimate}
Given the same conditions as the \Cref{cor:non-self-adjoint}, for $a:=h_\phi^{\frac \eta{p}} xh_\phi^{\frac {1-\eta}{p}} \in L^p(\M)$, $x\in\M^\sigma_a$,
\[
\|a-E_p^\phi(a)\|_{L^p(\M)} \leq 2p \sqrt{\frac{2}{\alpha}} \|x\|_{\Lip_\Gamma},
\]
where the Lipschitz seminorm \cite{gao_fisher_2020} is $\|x\|_{\Lip_\Gamma}
:= \max\{\|\Gamma(x,x)\|_\infty^{1/2}, \|\Gamma(x^*,x^*)\|_\infty^{1/2}\}$.
\end{cor}
\begin{proof}
To use the derivation trick, we first consider $\eta=\frac12$ and
self-adjoint $x=x^*$. The case of general $\eta\in[0,1]$ and
$x\in\M^\sigma_a$ follows by similar decomposition into self-adjoint and skew-adjoint parts as
in \Cref{cor:non-self-adjoint-ppp}, with an extra factor $2$.
On the slice $(\M_k,\tau_k)$, let
\[
  L^k_\varepsilon:=\frac{I-T^k_\varepsilon}{\varepsilon}
\]
be the $\varepsilon$-regularized generator. Its spectral gap is $\alpha_\varepsilon=\frac{1-e^{-\alpha\varepsilon}}{\varepsilon}$
by \Cref{rmk:t-reg-sg}. The $\Poin(p,p)$ inequality takes the form
\[
\|a_k- E_k(a_k)\|_{L^p(\M_k)}
\leq
\frac{p}{\sqrt{2\alpha_\varepsilon}}
\|a_k\|_{\Gamma_{L^k_{\varepsilon},p}}.
\]
The right-hand side can be expressed in terms of the derivation
$\partial^k_{\varepsilon}$ constructed in \Cref{cor:p-derivation}. In
Haagerup's $L^p$ notation, writing
$\Gamma_p^k=\Gamma_{L^k_{\varepsilon},p}$, we have
\begin{align*}
\sqrt{2}
\|a_k\|_{\Gamma_p^k}
&=
\Bigl\|
E_{\M_k}\bigl(
(\partial^k_{\varepsilon}(a_k))^*
\partial^k_{\varepsilon}(a_k)
\bigr)
\Bigr\|_{L^{p/2}(\M_k)}^{1/2}                              \\
&\le
\|\partial^k_{\varepsilon}(a_k)\|_{L^p(\mathscr M_k)}              \\
&=
\bigl\|
d_k^{1/2p}
\partial^k_{\varepsilon}(x_k)
d_k^{1/2p}
\bigr\|_{L^p(\mathscr M_k)}                                        \\
&\le
\|d_k^{1/2p}\|_{L^{2p}(\mathscr M_k)}^2
\|\partial^k_{\varepsilon}(x_k)\|_{\infty}                       =
\|\partial^k_{\varepsilon}(x_k)\|_{\infty}                         \\
&\le
2\sqrt{2}
\|\Gamma^k_{L_\varepsilon}(x_k,x_k)\|_{\infty}^{1/2}.
\end{align*}
The first equality follows from \Cref{cor:p-derivation}. The next
inequality uses the contractivity of conditional expectations. The third
line uses \cref{eq:p-derivation}. The fourth line follows from
H\"older's inequality and $\tau_k(d_k)=1$ since $d_k$ is density of a normal faithful state $\phi_k$. The last inequality follows from
\cref{eq:infty-estimate}.

Since the estimate is uniform in $k\ge1$, we argue as in the proof of
\Cref{thm:GNS-Lp-ppp} and take $\epsilon\downarrow 0$. Passing to the Haagerup reduction limit
$k\to\infty$ yields the desired estimate.
\end{proof}
\begin{remark}
\Cref{cor:lip-norm-estimate} may be compared with \cite[Theorem~5.12]{Gao_Junge_LaRacuente_Li_2025}, where a $\Poin(p,\infty)$ inequality with constant $O(\sqrt p)$ is obtained under the assumption $\alpha_{\mathrm{MLSI}}>0$, rather than an $L^2$ spectral gap.
In view of the discussion in the introduction, the two results should be regarded as independent.
\end{remark}
We follow the procedure of \cite{Gao_Junge_LaRacuente_Li_2025} for the sub-exponential concentration property.
\begin{definition}
For $x\in\M$, we write
\[
\bP_{\phi}(|x|>t)\kl\epsilon
\]
if there exists a projection $e\in\M$ such that
\[
\|exe\|_{\infty}\le t
\quad\text{and}\quad
\phi(1-e)\le\epsilon .
\]
\end{definition}

\begin{lemma}[\text{\cite[Lemma~5.14]{Gao_Junge_LaRacuente_Li_2025}}]\label{lem:cheb}
Let $x\in L_p(\M,\phi)$ with $1<p<\infty$. Then
\[
\bP_{\phi}(|x|>t)
\kl
2\Bigl(\frac{t}{4}\Bigr)^{-p}\|x\|_{L^p(\M)}^{p}\pl .
\]
\end{lemma}

\begin{cor}[Sub-exponential Concentration]\label[corollary]{cor:concentration}
Let $T_t=e^{-tL}$ be a GNS--$\phi$-detailed balance QMS with spectral gap $\alpha>0$ and $\phi$-preserving conditional expectation $E^\phi$ onto its fixed-point algebra. Then for any $x\in\M$ and large $t>0$,
\[
\bP_{\phi}(|x-E^\phi(x)|>t)
\le
2\exp\Bigl(-\frac{\sqrt{\alpha} t}{8\sqrt{2}e \|x\|_{\Lip_\Gamma}}\Bigr).
\]
\end{cor}

\begin{proof}
By \Cref{lem:cheb} and \Cref{cor:lip-norm-estimate},
\[
\bP_{\phi}(|x-E(x)|>t)
\le
2\Bigl(\frac{t}{4}\Bigr)^{-p}\|x-E(x)\|_{L_p(\M)}^{p}
\le
2\Bigl(\frac{8\sqrt{2}p\|x\|_{\Lip_\Gamma}}{\sqrt{\alpha} t}\Bigr)^{p}\pl .
\]
Optimizing in $p$ gives $p=\frac{\sqrt{\alpha} t}{8\sqrt{2}e\|x\|_{\Lip_\Gamma}}\ge3$ for $t$ large, which implies the claim.
\end{proof}
This is a generalization of the sub-exponential concentration obtained in the semi-commutative case \cite[Thm~2.7]{huang2021poincare}. Under only spectral gaps, one cannot expect tail decay faster than exponential, since this hypothesis does not imply stronger concentration phenomena. In contrast, \cite[Cor~5.15]{junge2015noncommutative} proved Gaussian concentration under the additional assumption of MLSI, which is strictly stronger.
\subsubsection{Birth--death process}
The following birth-death process shows that the $O(p)$ $\Poin(p,p)$ inequality in \Cref{main:GNS} is weaker than the $O(\sqrt p)$ Lipschitz-moment estimate under an MLSI assumption \cite[Theorem~5.12]{Gao_Junge_LaRacuente_Li_2025}.
We illustrate this by a violation of the geometric Talagrand inequality by considering a birth-death process with a constant spectral gap.
However, by \cite[Thm~6.10 and Rem.~6.9]{gao_fisher_2020}, $O(\sqrt p)$ Lipschitz-moment estimate implies a Wasserstein-$1$ bound and hence geometric Talagrand concentration. Therefore, the $O(\sqrt p)$ estimate cannot follow from a spectral gap alone, whereas the $O(p)$ estimate in \Cref{main:GNS} still holds.

\begin{definition}[Geometric Talagrand inequality]\label[definition]{def:talagrand}
Let $T_t=e^{-tL}$ be a QMS on $(\M,\tau)$ with fixed-point algebra $\N$ and conditional expectation $E\colon\M\to\N$.
We say that $T_t$ satisfies the geometric Talagrand inequality if there exists $c>0$ such that for every normal state $\rho$,
\begin{equation}\label{eq:tala}
\|\rho-E_*(\rho)\|_{\Gamma^*}
\le
c \sqrt{ D(\rho\|E_*(\rho))},
\end{equation}
where $\displaystyle\|\rho\|_{\Gamma^*}:=\sup\{\bigl|\tau(\rho x)\bigr|: E(x)=0, {x=x^*,\ \|x\|_{\Lip_\Gamma}\le1}\}$ and $E_*$ is the predual of $E$.
\end{definition}

\begin{example}[Birth-death Process]\label[example]{ex:bd}
  For a 1D chain of length $n$ and fixed $\beta>0$, we define the birth-death process edge-wise as follows
\begin{align*}
L_{j,j+1}(x)
=&e^{\beta/2}\bigl(\ket{j+1}\bra{j+1}x+x\ket{j+1}\bra{j+1}-2\ket{j+1}\bra{j}x\ket{j}\bra{j+1}\bigr)\\
 &+e^{-\beta/2}\bigl(\ket{j}\bra{j}x+x\ket{j}\bra{j}-2\ket{j}\bra{j+1}x\ket{j+1}\bra{j}\bigr).
\end{align*}
The total birth-death Lindbladian is
\[
L=\sum_{j=1}^{n-1} L_{j,j+1}.
\]
The Lindbladian $L$ satisfies GNS-$\mu$ detailed balance for the thermal state $\mu:=Z_\beta^{-1} \sum_{k=1}^n e^{-\beta k}\ket{k}\bra{k}$, where $Z_\beta:=\sum_{j=1}^n e^{-\beta j}$ is the partition function. Restricting to the diagonal part $L^{\mathrm{diag}}$ and fixed $\beta>0$, \cite{Miclo1999DiscreteHardy} gives the spectral gap $\alpha(L^\mathrm{diag})=O(1)$. We therefore establish $\PI(p,p)$ with constant $O(p)$ by \Cref{main:GNS}.

Nevertheless, the Talagrand's inequality \Cref{eq:tala} admits no constant
$c=c(\beta)$ independent of $n$.
Following the setup in \cite[Remark~6.9]{gao_fisher_2020}, consider $e_A$ a projection in region $A$ and density matrix $\rho_A:=\frac{e_A}{\Tr(e_A)}$. Choose $f=\sum_{j=1}^n \frac j{\sqrt{2\cosh{\beta/2}}} \ket j\bra j$ which gives  \[\|f\|_{\Lip_\Gamma}=\|\sum_j \bigl|[e^{\beta/4}(\ket{j}\bra{j+1}),f]\bigr|^2+\sum_j \bigl|[e^{-\beta/4}(\ket{j+1}\bra{j}),f]\bigr|^2\|^{1/2}
=1.\]
Let $A=\{1\}$ and $B=\{n\}$. Since $\rho_A$ and $\rho_B$ are supported on disjoint $A$ and $B$,
\[
\|\rho_A-\rho_B\|_{\Gamma^*}
=\sup_{x=x^*,\ \|x\|_{\Lip_\Gamma}\le1}\bigl|\Tr((\rho_A-\rho_B)x)\bigr|\ge\bigl|\Tr((\rho_A-\rho_B)f)\bigr|=\frac{n-1}{\sqrt{2\cosh{\beta/2}}}.
\]
Moreover, \[D\Big(\rho_{\{i\}}||E_*(\rho_{\{i\}})\Big)=D(\rho_{\{i\}}||\mu)=-\ln \mu ({\{i\}}),\quad   i\in\{1,...,n\}.\]
The right-hand side ${c(\beta)}(\sqrt{D(\rho_A||E_*(\rho_A))}+\sqrt{D(\rho_B||E_*(\rho_B))})=O(\sqrt{n})$, whereas the left-hand side $\|\rho_A-\rho_B\|_{\Gamma^*}$ of linear order in $n$.
This contradicts \Cref{eq:tala} for any fixed $c=c(\beta)$ as $n\to\infty$.
\end{example}
The birth-death process above provides an example showing that a subgaussian-type estimate, where $K$ is a universal constant
\[
\|f-E(f)\|_p \le Kc(\beta)^2\sqrt p \|f\|_{\Lip_\Gamma},
\]
fails uniformly in $n$ with a constant depending only on $\beta$.

This agrees with the result \cite[Thm~5.7]{Gao_Junge_LaRacuente_Li_2025} where they established the CMLSI $\alpha_c=O(1/n)$ replacing $c(\beta)^2$.
In contrast, $L^{\mathrm{diag}}$ has a spectral gap of $O(1)$, so the $\Poin(p,p)$ inequality with constant $O(p)$ still holds.

\subsection{Lipschitz diameter}
As an application of \Cref{main:GNS}, we obtain Lipschitz diameter and completely bounded (CB) Lipschitz diameter estimates for GNS-$\phi$-detailed-balanced QMSs in finite dimensions, depending on the spectral gap $\alpha$ and the spectral density of $\phi$.
\begin{cor}[Lipschitz diameter]\label[corollary]{cor:lip-cb-diam}
Let $\bM_d$ be a finite-dimensional matrix algebra equipped with the unnormalised trace $\Tr$ and large dimension (e.g. $d> 20$),
and denote by $\|\cdot\|_{S_p}:= \bigl( \Tr |\cdot|^p \bigr)^{1/p}$ the Schatten--$p$ norm.
Let $(T_t)_{t\ge0}$ be a quantum Markov semigroup satisfying GNS detailed balance with respect to
$\phi(\cdot)=\Tr(D_\phi \cdot)$ and spectral gap $\alpha>0$.
Then for any $x\in\bM_d$,
\[
\|x-E(x)\|_{\infty}
\le
\mathrm{diam}_{\Gamma}
\|x\|_{\Lip_\Gamma},
\]
where $\mathrm{diam}_{\Gamma}\le 2e\log(\lambda_{\min}^{-1})\sqrt{\frac{2}{\alpha}}
$ and $\lambda_{\min}$ denotes the minimal eigenvalue of $D_\phi$.

Moreover, invoking Pisier's vector-valued Schatten norm \cite{Pisier_1996} at each matrix level, the corresponding completely bounded diameter estimate is
\[
  \mathrm{diam}_{\Gamma}^{\mathrm{cb}}
  :=
  \sup_{m\ge1}\mathrm{diam}_{\Gamma^{(m)}}
  \le
  4e\log(\lambda_{\min}^{-1})\sqrt{\frac{2}{\alpha}} .
\]
\end{cor}

\begin{proof}
We write $D$ as shorthand for $D_\phi$ and, without loss of generality,
$x\in\mathring{\bM_d}$, where $\mathring{\bM_d}=\ker(E)$. Since $x
=
D^{-\frac1{2p}}
 \bigl(D^{\frac1{2p}} x D^{\frac1{2p}}\bigr)
 D^{-\frac1{2p}},$
the Hölder inequality and \Cref{cor:lip-norm-estimate} gives
\begin{align*}
\|x\|_{\infty}
&\le \|x\|_{S_p} \le
\bigl\|D^{-\tfrac1{2p}}\bigr\|_{\infty}^{2}
\bigl\|D^{\tfrac1{2p}} x D^{\tfrac1{2p}}\bigr\|_{S_p} \\
&\le (\lambda_{\min})^{-\tfrac1p}
\|\iota_{1/2}(x)\|_{p,\phi,1/2} \le
(\lambda_{\min})^{-\tfrac1p} 2p
\sqrt{\frac{2}{\alpha}} \|x\|_{\Lip_\Gamma},
\end{align*}
Now choose
$
p = \log\bigl(\lambda_{\min}^{-1}\bigr)\ge \log\bigl(d\bigr)>3
$,
so that $(\lambda_{\min})^{-1/p}=e$.  Substituting back leads to the desired inequality.

For the complete version, fix $m\ge1$ and consider the amplified
semigroup $\Id_m\otimes T_t$. Its fixed-point algebra is
$\bM_m\overline{\otimes}\N$, and its conditional expectation onto this
algebra is $\Id_m\otimes E$. Let
$\Gamma^{(m)}$ denote the gradient form of $\Id_m\otimes T_t$. Let $X\in\bM_m(\bM_d)$ be self-adjoint and satisfy
\[
  X=X^*\in\ker(\Id_m\otimes E)\subset\bM_m(\bM_d).
\]
For every $p\ge3$, Pisier's formula for self-adjoint elements
\cite{Pisier_1996} gives
\begin{equation}\label{eq:pisier-complete-diameter}
  \|X\|_{\bM_m(\bM_d)}
  =
  \sup_{\psi}
  \bigl\|(D_\psi^{1/(2p)}\otimes1)X
  (D_\psi^{1/(2p)}\otimes1)\bigr\|_{S_p^m[\bM_d]},
\end{equation}
where the supremum runs over the states $\psi$ on $\bM_m$ and $D_\psi$ denotes
the density of $\psi$ with respect to $\Tr_m$. i.e. $\psi(\cdot)=\Tr_m(D_\psi\cdot)$.

To estimate the RHS of
\cref{eq:pisier-complete-diameter}, we can expand $1=D^{\frac1{2p}}D^{-\frac1{2p}}$, and rewrite
\[(D_\psi^{1/(2p)}\otimes1)X
  (D_\psi^{1/(2p)}\otimes1)=(\Id_{S_p^m}\otimes\alpha_p)(Z_\psi),
\]
where $Z_\psi:=(D_\psi^{1/(2p)}\otimes D^{1/(2p)})X(D_\psi^{1/(2p)}\otimes D^{1/(2p)})$ and the map $\alpha_p$ is defined as
\[
  \alpha_p:S_p^d\longrightarrow\bM_d,
  \qquad
  \alpha_p(Y)=D^{-1/(2p)}YD^{-1/(2p)}.
\]
The completely bounded two-sided multiplication estimate for the canonical
operator-space structure on $S_p^d$ \cite{AST_1998__247__R1_0} yields
\begin{equation}\label{eq:alpha-p-cb}
  \|\alpha_p\|_{\mathrm{cb}(S_p^d,\bM_d)}
  \le
  \bigl\|D^{-1/(2p)}\bigr\|_{S_{2p}}^2
  =
  (\sum_j\lambda_j^{-1})^{1/p}.
\end{equation}
Using \cref{eq:alpha-p-cb} and the complete isometry
$S_p^m[S_p^d]=S_p^{md}$, we obtain estimate of RHS of \cref{eq:pisier-complete-diameter}:
\begin{equation}
\bigl\|(D_\psi^{1/(2p)}\otimes1)X
  (D_\psi^{1/(2p)}\otimes1)\bigr\|_{S_p^m[\bM_d]}
  \le
  (\sum_j\lambda_j^{-1})^{1/p}
  \|Z_\psi\|_{S_p^{md}}.
  \label{eq:complete-diameter-factorization}
\end{equation}

The amplified semigroup $\Id_m\otimes T_t$ has spectral gap $\alpha$ relative
to $\Id_m\otimes E$ by \Cref{prop:Gap-tensor}. Since $X=X^*$ and
$(\Id_m\otimes E)(X)=0$, the self-adjoint $\eta=\frac12$ estimate from the
proof of \Cref{cor:lip-norm-estimate}, applied to the state
$\psi\otimes\phi$, gives
\[
  \|Z_\psi\|_{S_p^{md}}
  \le
  p\sqrt{\frac{2}{\alpha}}
  \|X\|_{\Lip_{\Gamma^{(m)}}}.
\]
Combining this estimate with \cref{eq:pisier-complete-diameter,eq:complete-diameter-factorization}
and then taking the supremum over $\psi$ yields
\begin{equation}\label{eq:complete-diameter-p}
  \|X\|_{\bM_m(\bM_d)}
  \le
  p(\sum_j\lambda_j^{-1})^{1/p}
  \sqrt{\frac{2}{\alpha}}
  \|X\|_{\Lip_{\Gamma^{(m)}}}.
\end{equation}
Set $A=\sum_j\lambda_j^{-1}$ and choose $p=\log A$ to minimize the RHS of \cref{eq:complete-diameter-p}. We have
\[
  \mathrm{diam}_{\Gamma}^{\mathrm{cb}}
  \le
  e\log A\sqrt{\frac{2}{\alpha}}
  \le
  2e\log(\lambda_{\min}^{-1})\sqrt{\frac{2}{\alpha}},
\]
where we used
$A\le d\lambda_{\min}^{-1}\le\lambda_{\min}^{-2}$ in the last inequality.
Finally, for a general $X$, the decomposition into its self-adjoint and
skew-adjoint parts, as in \Cref{cor:non-self-adjoint-ppp}, contributes a factor
$2$. Taking the supremum over $m\ge1$ therefore yields the desired completely bounded estimate.
\end{proof}
Returning to the birth-death process in \Cref{ex:bd} for fixed $\beta>0$, we can apply the above corollary to obtain a diameter estimate for the semigroup. The minimal eigenvalue of the thermal state $\mu$ is $\lambda_{\min} = Z_\beta^{-1} e^{-\beta n}$, which gives
\[
  \operatorname{diam}_\Gamma=O(n),
  \qquad
  \operatorname{diam}_\Gamma^{\mathrm{cb}}=O(n).
\]

\appendix
\addtocontents{toc}{\protect\setcounter{tocdepth}{1}}
\crefalias{section}{appendix}

\section{Markov dilation and proof of \texorpdfstring{\Cref{prop:convex-chain}}{Prop}}\label{app:chain}
Before proving the noncommutative convex-chain estimate for gradient forms, we first
recall the framework of Markov dilation, which lifts a quantum Markov semigroup to a
multiplicative level.
This allows the Dirichlet form to be written as a square and enables the
noncommutative chain rule.

Markov dilation embeds quantum Markov semigroups into larger operator-algebraic structures that carry probabilistic interpretations. The foundational work of Kümmerer and Maassen \cite{Kümmerer_Maassen_1987} introduced noncommutative dilation in which the enlarged algebra remains essentially commutative, and the first author and Mei \cite{junge-mei} extended this framework to the setting of tracial von Neumann algebras. In \cite{wu2023}, the author provided a probabilistic construction of Markov dilation. In quantum physics, such dilations are interpreted as providing a reversible microscopic model of systems interacting with thermal environments. We recall the definition from \cite[Section~2.2]{junge-mei}.

\begin{definition}[Markov dilation]\label{def:dilation}
Let $(T_t)_{t\ge0}$ be a quantum Markov semigroup on a tracial von Neumann algebra $(\M,\tau)$.
We say that $(T_t)$ admits a Markov dilation if there exists an increasing family of tracial von Neumann algebras
\[
  \M \subset \wt{\M}_t \subset \wt{\M}, \quad t \ge 0,
\]
that is, an increasing filtration of $\wt{\M}_t$, together with
\begin{itemize}
  \item normal, trace-preserving conditional expectations $E_t:\wt{\M}\to\wt{\M}_t$,
  \item $*$--homomorphisms $\pi_t:\M \to \wt{\M}$ that are adapted to the filtration, i.e. $\pi_t(\M)\subset\wt\M_t$, for each $t\ge0$,
\end{itemize}
such that for all $0\le s<t<\infty$ and all $x\in\M$,
\begin{equation}\label{eq:dilation}
  E_s\bigl(\pi_t(x)\bigr) = \pi_s\bigl(T_{t-s}(x)\bigr).
\end{equation}
Equivalently, condition \eqref{eq:dilation} is expressed by the following commuting diagram.
\begin{equation}\label{eq:dilation-diagram}
  \begin{tikzcd}[row sep=1.5em, column sep=3em]
    \M \arrow[r, "\pi_t"] \arrow[d, swap, "T_{t-s}"]
      & \wt{\M} \arrow[d, "E_s"] \\
    \M \arrow[r, swap, "\pi_s"]
      & \wt{\M}_s
  \end{tikzcd}
\end{equation}
\end{definition}
We now prove \Cref{prop:convex-chain}, restated below for convenience.
\begin{proposition*}[Convex-chain estimate for gradient form]
Given the conditions of Section~\ref{sec:preliminaries} and lemma \ref{lem:klein},
for every self-adjoint element
$x\in\Dom(L)\subset L_{2}(\M,\tau)$,
\[
  \mathscr E\bigl(\varphi(x)\bigr)
   \le
  \tau   \bigl(\Gamma(x,x)  \psi(x)\bigr),
\]
where
$\mathscr E(f):=\tau   \bigl(\Gamma(f,f)\bigr)$ is the Dirichlet form
associated with $L$.
\end{proposition*}
\begin{proof}
Use Markov dilation from \Cref{def:dilation}, for $x\in\M$, and from the commuting diagram \eqref{eq:dilation-diagram} with $s=0$ $T_t=E_0\circ\pi_t|_{\M}$ and that $T_t$ is $\tau$-symmetric as in equation \eqref{eq:tau-sym},
\begin{align*}
\mathscr{E}(x)
      &:= \tau   \bigl(x^{*}  {L}(x)\bigr) = \lim_{t\to0}
         \tau   \Bigl(
           x^{*}  \frac{1-e^{-t{L}}}{t}  x
         \Bigr) = \lim_{t\to0}
         \frac{\tau(x^{*}x)-\tau(x^{*}T_t(x))}{t} \\
      &= \lim_{t\to0}
         \frac{1}{2t}
         \tau   \Bigl(
           x^*x+T_t(x^*x)-T_t(x)^*  x-x^*  T_t(x)
         \Bigr) \\
    &=\lim_{t\to0}
         \frac{1}{2t}
         \tau   \Bigl(
           x^*x+\pi_t(x^*x)-\pi_t(x)^*  x-x^*  \pi_t(x)
         \Bigr) \\
      &= \lim_{t\to0}
         \frac{1}{2t}
         \tau   \Bigl(\bigl(x-\pi_t(x)\bigr)^*\bigl(x-\pi_t(x)\bigr)\Bigr).
\end{align*}
Here the first line is the usual definition of the generator equation \eqref{eq:lindbladian}.
The second line uses the $\tau$-symmetry of the semigroup, namely $\tau(x^{*}T_t(x))=\tau(T_t(x)^{*}x)$ and $\tau(T_t(x^*x))=\tau(T_t(1)^{*}x^*x)=\tau(x^*x)$.
In the third line we replace $T_t$ by $E_0\pi_t$ and drop $E_0$ inside the trace because $\tau\circ E_0=\tau$ and $x\in\M=\wt{\M}_0$ so $E_0(x)=x$.
Finally, since $\pi_t$ is a $*$ homomorphism, we reorganize the expression into the squared difference in the last line.
Now restrict to $x\in\M_{sa}$, so that $\varphi(x),\psi(x)\in\M_{sa}$. We then apply the same argument to $\varphi(x)$ and use \Cref{lem:klein} to obtain
\[
\begin{aligned}
\mathscr{E}   \bigl(\varphi(x)\bigr)
  &= \lim_{t\to0}\frac{1}{2t}
     \tau   \Bigl(
       \bigl(\varphi(x)-\pi_t   \bigl(\varphi(x)\bigr)\bigr)^2
     \Bigr)                            \\
  &= \lim_{t\to0}\frac{1}{2t}
     \tau   \Bigl(
       \bigl(\varphi(x)-\varphi   \bigl(\pi_t(x)\bigr)\bigr)^2
     \Bigr)                            \\
  &\le
     \lim_{t\to0}\frac{1}{4t}
     \tau   \Bigl(
       \bigl(x-\pi_t(x)\bigr)^{2}
       \bigl(\psi(x)+\psi   \bigl(\pi_t(x)\bigr)\bigr)
     \Bigr)                            \\
  &= \lim_{t\to0}\frac{1}{4t}
     \tau   \Bigl(
       x^{2}\psi(x)-x  \pi_t(x)\psi(x)-\pi_t(x)  x\psi(x)+\pi_t(x)^{2}\psi(x)
     \Bigr)                            \\
  &\quad+
     \lim_{t\to0}\frac{1}{4t}
     \tau   \Bigl(
       x^{2}\pi_t   \bigl(\psi(x)\bigr)
       -x  \pi_t(x)\pi_t   \bigl(\psi(x)\bigr)
       -\pi_t(x)  x\pi_t   \bigl(\psi(x)\bigr)
       +\pi_t(x)^{2}\pi_t   \bigl(\psi(x)\bigr)
     \Bigr)                            \\
    &=  \lim_{t\to0}\frac{1}{4t}
     \tau   \Bigl(
       x^{2}\psi(x)-T_t(x)\psi(x)x-T_t(x)  x\psi(x)+T_t(x^{2})\psi(x)
     \Bigr)                            \\
  &\quad+
     \lim_{t\to0}\frac{1}{4t}
     \tau   \Bigl(
       x^{2}\pi_t   \bigl(\psi(x)\bigr)
       -x  \pi_t(x\psi(x))
       -x\pi_t   \bigl(\psi(x)x\bigr)
       +\pi_t   \bigl(x^{2}\psi(x)\bigr)
     \Bigr)                            \\
         &=  \lim_{t\to0}\frac{1}{4t}
     \tau   \Bigl(\bigl(
       x^{2}-xT_t(x)-T_t(x)  x+T_t(x^{2})\bigr)\psi(x)
     \Bigr)                            \\
     &\quad+
     \lim_{t\to0}\frac{1}{4t}
     \tau   \Bigl(
       T_t   \bigl(x^2\bigr)\psi(x)
       -T_t(x)  x\psi(x)
       -xT_t(x)\psi(x)
       +x^{2}\psi(x)
     \Bigr)                            \\
  &= \tau   \bigl(\Gamma(x,x)  \psi(x)\bigr),
\end{aligned}
\]
\noindent
In these manipulations we repeatedly use the cyclicity of the trace, the multiplicativity and trace preservation of $\pi_t$, the $\tau$-symmetry of $T_t$, and
$\tau\circ E_0=\tau$.  After canceling the mixed terms with those
properties, the remaining $t$--difference quotient is recognized as the
gradient form $\Gamma(x,x)$ defined in~\eqref{eq:Gamma}.
\noindent
Notice that we \emph{never} assumed $\pi_t(x)$ and $x$ to be exchangeable;
$\pi_t$ is not required to be an automorphism on~$\M$.
The argument relies only on the elementary properties already listed
(multiplicativity of $\pi_t$
and trace preservation), and not on any further commutation relations.
The price for working under these minimal hypotheses is the
lengthy algebraic expansion needed to isolate the gradient form term.
\end{proof}
\section{Spectral gaps in AFP}
This section identifies the spectral gap for tensor, direct-sum, and amalgamated product quantum Markov semigroups.
\begin{proposition}[Spectral gap for $T^1_t\otimes T^2_t$]\label{prop:Gap-tensor}
Let $(T^i_t)_{t\ge0}$ be quantum Markov semigroups on $\M_i$ with conditional expectations $E_i:\M_i\to\N_i$ onto the fixed-point algebras $\N_i$ with spectral gaps $\alpha_i>0$, i.e.
\[
\|T^i_t(x)-E_i(x)\|_{2}
 \le e^{-\alpha_i t}  \|x-E_i(x)\|_{2},
\quad   x\in\M_i, t\ge0.
\]
On
$\M:=\M_1\otimes\M_2$ consider the product semigroup $T_t
:= T^1_t\otimes T^2_t,$
with conditional expectation $E:=E_1\otimes E_2:\M_1\otimes\M_2\to\N_1\otimes\N_2$.  Then
\[
\|T_t(x)-E(x)\|_{2}
 \le
e^{-\min\{\alpha_1,\alpha_2\}  t}
\|x - E(x)\|_{2},
\quad   x\in\M, t\ge0,
\]
i.e. the spectral gap of $T^1_t\otimes T^2_t$ is $\min\{\alpha_1,\alpha_2\}$.
\end{proposition}
\begin{proof}
Regard $E_i$ as the orthogonal projection from
$\H_i:=L^2(\M_i,\tau_i)$ onto $L^2(\N_i,\tau_i)$, and set
$\mathring\H_i:=\ker E_i$.  Since
$T_t^iE_i=E_iT_t^i=E_i$, the orthogonal decomposition
\[
\ker(E_1\otimes E_2)
=
(\mathring\H_1\otimes\H_2)
\oplus
(\Ran E_1\otimes\mathring\H_2)
\]
reduces $T_t-E$.  On its two summands, the restrictions of $T_t-E$ are,
respectively,
\[
(T_t^1-E_1)\otimes T_t^2
\quad\text{and}\quad
E_1\otimes(T_t^2-E_2).
\]
It follows from the $L^2$-contractivity of $T_t^1$ and $T_t^2$ that
\begin{align*}
\|T_t-E\|_{B(\ker E)}
&\le
\max\bigl\{\|T_t^1-E_1\|,\|T_t^2-E_2\|\bigr\} \\
&\le
e^{-\min\{\alpha_1,\alpha_2\}t}.
\end{align*}
This proves the asserted estimate.

Finally, the restrictions of $T_t-E$ to
$\mathring\H_1\otimes\Ran E_2$ and
$\Ran E_1\otimes\mathring\H_2$ coincide with $T_t^1-E_1$ and
$T_t^2-E_2$, respectively.  Hence the spectral gap of the product semigroup
cannot exceed either $\alpha_1$ or $\alpha_2$.  The preceding estimate therefore
shows that it is exactly $\min\{\alpha_1,\alpha_2\}$.
\end{proof}

\begin{proposition}[Spectral gap in $\M_1\oplus\M_2$ via semigroups]\label{prop:Gap-sum}
Under the setup of \Cref{prop:Gap-tensor}, on
$\M:=\M_1\oplus\M_2$ define $T_t:=T^1_t\oplus T^2_t$ and $E:=E_1\oplus E_2$.  Then
\[
\|T_t(x)-E(x)\|_2
 \le
e^{-\min\{\alpha_1,\alpha_2\}t}
  \|x-E(x)\|_2 \qquad
   x\in\M_i, t\ge0,
\]
i.e. the spectral gap is $\min\{\alpha_1,\alpha_2\}$.
\end{proposition}
\begin{proof}
For $x=(x_1,x_2)$ in $\H$, $T_t(x)-E(x)
=\bigl(T^1_t(x_1)-E_1(x_1),  T^2_t(x_2)-E_2(x_2)\bigr),$
so we have
\begin{align*}
  \|T_t(x)-E(x)\|_2^2
&=\sum_{i=1}^2\|T^i_t(x_i)-E_i(x_i)\|_2^2\\
&\le
\sum_{i=1}^2e^{-2\alpha_i t}  \|x_i-E_i(x_i)\|_2^2
\le
e^{-2\min\{\alpha_1,\alpha_2\}t}  \|x-E(x)\|_2^2.
\end{align*}
Taking square roots gives the desired estimate.
\end{proof}

\begin{cor}[Spectral gap in $\M_1*_\N\M_2$ via semigroups]\label{cor:gap}
Under the setup of \Cref{prop:Gap-tensor},
let $\M:=\M_1*_\N\M_2$ carry the free-product semigroup $T_t$ \cite{boca_completely_1993} with joint expectation $E$.  Then for all $x\in L_2(\M)$,
\[
\|T_t(x)-E(x)\|_{2}
 \le
e^{-\min\{\alpha_1,\alpha_2\}t}  \|x-E(x)\|_{2}, \quad   x\in\M,t\ge0
\]
and hence the spectral gap is $\min\{\alpha_1,\alpha_2\}$.
\end{cor}
\begin{proof}
Recall the orthogonal decomposition of the GNS-space \eqref{eq:hilbert-free}
\[
L_2(\M)
 =
L_2(\N)
 \oplus
\bigoplus_{\substack{w=i_1\cdots i_k\\ i_j\in\{1,2\},  i_j\neq i_{j+1}}}
\H_w,
\]
where each $\H_w$ is the closed linear span of reduced words $x_{i_1}\cdots x_{i_k}$ with $E_{i_j}(x_{i_j})=0$.  On each summand $\H_w$, $T_t$ acts by
\[
T_t(x_{i_1}\cdots x_{i_k})
= \pi_{i_1}\circ T^{(i_1)}_t(x_{i_1}) \cdots \pi_{i_k}\circ T^{(i_k)}_t(x_{i_k}),
\]
For simplicity, we give the proof in the scalar free-product case
$\N=\mathbb C1$.
For general $\N$, the same proof applies to the AFP decomposition by replacing ordinary tensor products with relative tensor products $\overline{\otimes}_{\N}$ and using $\N$-bimodularity together with the induced tensor-operator norm estimate.

By orthogonality and the tensor-product estimate (\Cref{prop:Gap-tensor}),
\[
\|T_t(x_{i_1}\cdots x_{i_k})\|_2
 =
\prod_{j=1}^k\|T^{(i_j)}_t(x_{i_j})\|_2
 \le
\prod_{j=1}^k\bigl(e^{-\alpha_{i_j}t}\|x_{i_j}\|_2\bigr)
 =
e^{-\bigl(\sum_j\alpha_{i_j}\bigr)t}
\|x_{i_1}\cdots x_{i_k}\|_2.
\]
In particular, the slowest decay on a word arises from its last letter, so
\[
\|T_t|_{\H_w}\|
 \le
e^{-\min\{\alpha_{i_k}\}  t}
 \le
e^{-\min\{\alpha_1,\alpha_2\}  t}.
\]
Since $L_2(\mathring{\M})=L_2(\M)\ominus L_2(\N)$ is the orthogonal direct sum of all the $\H_w$, we conclude
\[
\|T_t - E\|_{B(L_2(\M)\ominus L_2(\N))}
 =
\sup_w\|T_t|_{\H_w}\|
 \le
e^{-\min\{\alpha_1,\alpha_2\}t}.
\]
Equivalently,
\[
\|T_t(x)-E(x)\|_2
 \le
e^{-\min\{\alpha_1,\alpha_2\}t}  \|x-E(x)\|_2,
\]
and by definition this decay rate and Prop.~\ref{prop:Gap-sum} is exactly the spectral gap $\min\{\alpha_1,\alpha_2\}$.
\end{proof}

\section{Properties of \texorpdfstring{$\epsilon$}{epsilon}-regularized Lindbladians}\label{app:regular}
\subsection{Spectral gaps}
Under the notation of the main text, let $L$ be the generator of the
semigroup $(T_t)_{t\ge0}$, satisfying detailed balance (either
GNS or with respect to the trace), and let
$\alpha>0$ denote its spectral gap. We verify that if $x\in \Dom(L)$, that the $\epsilon$-regularization

$$\displaystyle L_\epsilon(x) := \frac{L}{1+\epsilon L}(x)\xrightarrow[{\epsilon\downarrow0}]{} L(x)$$
in the strong* operator topology.
We also verify that $L_\epsilon$
inherit the detailed balance property of $L$ and satisfy the lower bound
$\alpha_\epsilon \ge \alpha(1-\epsilon\alpha)$ for their spectral gaps.

The preservation of detailed balance follows directly from the integral
representation \cite[Prop.~2.5]{cipriani_derivations_2003},
\[
L_\epsilon
= \frac{1}{\epsilon}(I-R_\epsilon)
= \frac{1}{\epsilon}\bigl(I-\int_0^\infty e^{-t}T_{\epsilon t}dt\bigr),
\]
since $(T_t)_{t\ge0}$ satisfy detailed balance.

For the spectral gap, we restrict to $(\M_k,\tau_k)$ with $k\ge1$, where the
$\PI(p,p)$ inequality holds. Since $L_{k,\epsilon}$ is
$\tau_k$-symmetric, we apply spectral theory. If $\text{Spec}\bigl(L_k|_{\Ran L_k}\bigr)\subset[\alpha,\infty)$,
then by functional calculus
\[
\text{Spec}\bigl(L_{k,\epsilon}|_{\Ran L_{k,\epsilon}}\bigr)
= \Bigl\{ \tfrac{\lambda}{1+\epsilon\lambda} : \lambda\in\text{Spec}(L_k) \Bigr\}
\subset \Bigl[\tfrac{\alpha}{1+\epsilon\alpha},\infty\Bigr).
\]
Hence, the spectral gap gives
\[
\alpha_\epsilon = \frac{\alpha}{1+\epsilon\alpha} \ge \alpha(1-\epsilon\alpha),
\]
where the last inequality follows from $\tfrac{1}{1+x}\ge 1-x$ for $x\ge0$.

\begin{remark}\label{rmk:t-reg-sg}
    We can also repeat the spectral gaps for the $\varepsilon$-regularization in \cite{JungeRicardShlyakhtenko2025}, where
    \[
    L_\varepsilon(x):=\frac{1-T_\varepsilon(x)}{\varepsilon}\xrightarrow[{\varepsilon\downarrow0}]{} L(x),
    \]
in the strong* operator topology for $x\in \Dom(L)$. Assuming $L$ have spectral gap $\alpha>0$, the corresponding spectral gap gives
\[
\alpha_\varepsilon=\frac{1-e^{-\alpha \varepsilon}}{\varepsilon}\ge \frac{\alpha}{1+\varepsilon\alpha},
\]
which converges to $\alpha$ as $\varepsilon\downarrow 0$.
\end{remark}
Moreover, for a convex family $\overline L=\sum_j \lambda_j L_{\epsilon_j}$ with
$\epsilon_j>0$, $\lambda_j\ge0$, and $\sum_j\lambda_j=1$, detailed balance holds
since convex combinations of detailed-balance generators $L_{\epsilon_j}$ remain self-adjoint with respect to the relevant inner product.
By functional calculus, the spectral gap is
\[
\overline\alpha_\epsilon = \sum_j \lambda_j \frac{\alpha}{1+\epsilon_j\alpha},
\]
and the inequality
\[
\overline\alpha_\epsilon \ge \alpha \sum_j \lambda_j (1-\epsilon_j\alpha).
\]
\subsection{Strong convergence on the generator domain}
It is standard to check strong convergence in the domain and thus $p$-norm convergence. Indeed,
$p$-norm convergence follows from $\sigma$-strong convergence for bounded sequences
(or from strong convergence for bounded nets), as stated in the next lemma.
\begin{lemma}\label{lem:p-norm-convergence}
Let $(\M,\phi)$ be a von Neumann algebra with normal faithful state $\phi$ and
let $(x_\alpha) \subset \M$ be a bounded net (or sequence) converging strongly
(resp., $\sigma$-strongly) to $x \in \M$. Then for every $1 \leq p < \infty$ and every $\eta \in [0,1]$, in
Kosaki's notation $\iota_\eta(x_\alpha)\to \iota_\eta(x)$ in $L_p(\M,\phi)_\eta$-norm.
Moreover, if $\iota_\eta(x_\alpha)\to \iota_\eta(x)$ in $L_p(\M,\phi)_\eta$-norm for $1\le p<\infty$, then it is true for all $0<p<\infty$.

\end{lemma}
\begin{proof}
This follows directly from the Riesz-Thorin interpolation theorem, as shown in \cite[Lemma~2.3]{Junge2002Doobs}.
\end{proof}
\begin{proposition}\label{app:dom-strong}
Let $L\ge 0$ be a self-adjoint operator on $L^2(\M)$,
then $L_\epsilon (x)\to L(x)$ in the strong$^*$ topology for every
$x\in\Dom(L)$. Thus in Kosaki's notation, $\iota_\eta(L_\epsilon(x))\to \iota_\eta(L(x))$ in $L_p(\M,\phi)_\eta$-norm for $x\in\Dom(L)$.
Moreover, assuming that $\A= \Dom(L)\cap \M$ is a weak$^*$-dense $*$-algebra, $\Gamma_\epsilon(x,y)\to \Gamma(x,y)$ strongly for all $x,y\in\A$. Similarly, in Kosaki's notation, it converges in the $p$-norm.
\end{proposition}
\begin{proof}
Since $L\ge 0$ is self-adjoint, its spectral resolution satisfies
\[
L = \int_{0}^{\infty} t  dE_t .
\]
For $x\in\Dom(L)$, the functional calculus gives
\[
L_\epsilon x
= \int_{0}^{\infty} \frac{t}{1+\epsilon t}  dE_t x,
\qquad
L x
= \int_{0}^{\infty} t  dE_t x.
\]
Hence
\[
\|L_\epsilon x - Lx\|_{L^2(\M)}^{2}
= \int_{0}^{\infty}
   \Big|\frac{t}{1+\epsilon t} - t \Big|^{2}  d\|E_t x\|^{2}
= \int_{0}^{\infty}
   \frac{\epsilon^{2} t^{4}}{(1+\epsilon t)^{2}}  d\|E_t x\|^{2}.
\]
The integrand converges pointwise to $0$ as $\epsilon\downarrow 0$, and for
all $t\ge 0$, we have the bound
\[
0 \le \frac{\epsilon^{2} t^{4}}{(1+\epsilon t)^{2}} \le t^{2}.
\]
Since $x\in\Dom(L)$ implies $\int_{0}^{\infty} t^{2}  d\|E_t x\|^{2} < \infty,
$
the Dominated Convergence Theorem gives
$\|L_\epsilon x - Lx\|_{L^2(\M)}\to 0$.
Moreover, $\|L_\epsilon x\|_\infty\le\|Lx\|_\infty$,
because $(I+\epsilon L)^{-1}$ is a contraction on $\M$.
Uniform boundedness together with $L^2$-convergence therefore yields
strong convergence.
Thus $L_\epsilon \to L$ strongly on $\Dom(L)$.
It follows from $L_\epsilon$ being self-adjoint that
$L_\epsilon^* (x) \to L^* (x)$ for every $x\in\Dom(L)$. Hence, the convergence
is strong$^*$ on $\Dom(L)$.

For $x,y\in\A$, the above convergence permits a direct comparison of $ \Gamma_\epsilon(x,y)$ and $ \Gamma(x,y)$.

Given $  L_\epsilon(x)\to   L(x)$ strongly and $  L_\epsilon(x)^*\to   L(x)^*$ strongly, we obtain
\[
\|  L_\epsilon(x)^* y -   L(x)^* y\|_{L^2( {\M})}\to 0,
\qquad
\|x^*   L_\epsilon(y) - x^*   L(y)\|_{L^2( {\M})}\to 0 .
\]
Since $\A$ is a $*$-algebra, $x^* y \in \A \subseteq \Dom(  L)$, so $  L(x^*y)$ is well defined. By the definition of strong$^*$ convergence,
\[
\|  L_\epsilon(x^* y) -   L(x^* y)\|_{L^2( {\M})}\to 0 .
\]
All terms appearing in \cref{eq:regulated-Gamma} converge strongly. Hence $ \Gamma_\epsilon(x,y)\to  \Gamma(x,y)$ strongly for all $x,y\in\A$.

We also check that $\Gamma_\epsilon(x,x)$ is uniformly bounded in $\M$ as by \cite[Prop~2.5]{cipriani_derivations_2003}: $(1+\epsilon L)^{-1}$ is bounded, completely positive normal contractions on $\M$, by interpolation we have $\|L_{\epsilon,p}(x)\|_{L^p(\M)}\le \|L_p(x)\|_{L^p(\M)}$. For $x\in\Dom(L)$ we have $L(x)\in\M$, and therefore each term appearing in
\cref{eq:regulated-Gamma} is uniformly bounded in $L^{p/2}(\M)$.
By \Cref{lem:p-norm-convergence}, we obtain the $p$-norm convergence stated for $1<p<\infty$.
\end{proof}
\begin{remark}\label{rmk:unif-bdd}
From $L^1$ convergence as in \cref{eq:gamma-l1-convergence} and uniform boundedness of $\Gamma_\epsilon(x,x)$ in $\M$ for $x\in\Dom(L)$, we can apply the interpolation trick in \cite{Junge2002Doobs} to check $L^{p}$ norm convergence for $1\le p<\infty$.
\end{remark}

\section{General case: noncommutative diffusion semigroup}\label{app:conv}
We address the complication that $\Gamma(x,y)$ involves only
weak$^*$ convergence. To obtain convergence in the $L^p$-norm, one cannot pass
directly to the limit $\epsilon\downarrow0$. By an argument due to the unpublished work of the first author and collaborators \cite{JungeRicardShlyakhtenko2025}, one can form convex combinations
of the regularized terms $\sum_j\lambda_j\Gamma_{\epsilon_j}(x,y)$ that converge in $\|\cdot\|_{L^1(\M)}$
to $\Gamma(x,y)$, for $\epsilon_j>0,\lambda_j\ge0$ and $\sum_j\lambda_j=1$. In this way, norm convergence is recovered at the expense
of convexification.
We begin by introducing some definitions and results from \cite{JungeRicardShlyakhtenko2025}.
\begin{definition}[noncommutative Diffusion QMS]
A semigroup $(T_t)$ is called a noncommutative diffusion process if for every $x\in\Dom(L)\cap \M$ the gradient form satisfies $\Gamma(x,x)\in L^1(\M)$.
\end{definition}
\begin{proposition}
    $(T_t)$ is noncommutative diffusion if
    \begin{equation}\label{eq:gamma-l1-convergence}
        \Gamma_\epsilon(x,x)\to \Gamma(x,x)\in L^1(\M)
    \end{equation}
 in $\sigma(L_1(\M),\M)$ topology. This implies $\Gamma(x,x)\in \overline{\mathrm{conv}(\Gamma_\epsilon)}^{\|\|_{L^1(\M)}}$ by Hahn-Banach. Thus, there exists a family of Lindbladians $L_{\epsilon_j}$, and the gradient form of convex combinations of Lindbladians also converges
    \[\Gamma_{\sum_j\lambda_j L_{\epsilon_j}}(x,y)=\sum_j\lambda_j\Gamma_{\epsilon_j}(x,y)\to\Gamma(x,y)\]
in $\|\cdot\|_{L^1(\M)}$ as $\{\epsilon_j\}\downarrow0$
, for some $\epsilon_j>0,\lambda_j\ge0$ and $\sum_j\lambda_j=1$.
\end{proposition}
We now complete the proof of the general case. We first convexify the regularized
generators, apply the Poincar\'e inequality on each finite slice, and
then take the limits in the order $k\to\infty$ followed by
$n\to\infty$.
Writing $q=p/2$, for the fixed pair $(y,y)$, choose finite convex
combinations
\[
\overline L_n
:=
\sum_{j=1}^{m_n}\lambda_{n,j}\widehat L_{\epsilon_{n,j}},
\qquad
C_n
:=
\Gamma_{\overline L_n}
=
\sum_{j=1}^{m_n}\lambda_{n,j}
\widehat\Gamma_{\epsilon_{n,j}},
\]
where $\lambda_{n,j}\ge0$, $\sum_j\lambda_{n,j}=1$, and
$\delta_n:=\max_j\epsilon_{n,j}\to0$, such that
\begin{equation}\label{eq:Lp-conv-convex}
\left\|
\iota_{\frac12}\bigl(C_n(y,y)-\widehat\Gamma(y,y)\bigr)
\right\|_{q,\widehat\phi,\frac12}
\longrightarrow0.
\end{equation}
The preceding Hahn--Banach convexification gives $L^1$-convergence, and
the interpolation argument in \Cref{rmk:unif-bdd} gives the $L^q$-convergence in \cref{eq:Lp-conv-convex}.  The spectral gap of $\overline L_n$ is
\[
\overline\alpha_n
:=
\sum_{j=1}^{m_n}\lambda_{n,j}
\frac{\alpha}{1+\epsilon_{n,j}\alpha},
\]
and $\overline\alpha_n\to\alpha$ because $\delta_n\to0$.

For each $k$, let
$\overline L_{n,k}:=\overline L_n|_{\M_k}$.  Each regularized generator
commutes with $\cE_k$, so $\overline L_{n,k}$ is a
$\tau_k$-symmetric Lindbladian with spectral gap
$\overline\alpha_n$ and the same fixed-point expectation $E_k$ as
$L_k$.  Applying the tracial $\Poin(p,p)$ inequality to
$\overline L_{n,k}$ gives
\begin{align*}
A_k
&:=
\left\|\iota_{\frac12}\bigl(x_k-E_k(x_k)\bigr)
\right\|_{p,\phi_k,\frac12}
\\
&\le
\frac{p}{\sqrt{2\overline\alpha_n}}
\left\|
\iota_{\frac12}\Bigl(
\Gamma_{\overline L_{n,k},\frac12}^{(p,k)}(x_k,x_k)
\Bigr)
\right\|_{q,\phi_k,\frac12}^{1/2}
\\
&=
\frac{p}{\sqrt{2\overline\alpha_n}}
\left\|
\cE_{k,1}\!\left(
\iota_{\frac12}\bigl(C_n(y,y_k)\bigr)
\right)
\right\|_{q,\widehat\phi,\frac12}^{1/2}
\\
&\le
\frac{p}{\sqrt{2\overline\alpha_n}}
\left\|
\iota_{\frac12}\bigl(C_n(y,y_k)\bigr)
\right\|_{q,\widehat\phi,\frac12}^{1/2}.
\end{align*}
Here the equality follows from the bimodularity of $\cE_k$ and its
intertwining with each $\widehat L_{\epsilon_{n,j}}$, while the last
inequality follows from the contractivity of its $L^q$ extension.

For fixed $n$, since $C_n$ is a finite convex combination of bounded
regularized gradient forms and $y_k\to y$ $\sigma$-strongly,
\[
\iota_{\frac12}\bigl(C_n(y,y_k)\bigr)
\xrightarrow{k\to\infty}
\iota_{\frac12}\bigl(C_n(y,y)\bigr)
\qquad\text{in }L^q(\widehat\M).
\]
On the other hand, \cref{eq:1,eq:2} gives convergence of the left-hand side of the Poincar\'e inequality.
Therefore, taking $k\to\infty$ with $n$ fixed yields
\[
\left\|\iota_{\frac12}\bigl(x-E(x)\bigr)
\right\|_{p,\phi,\frac12}
\le
\frac{p}{\sqrt{2\overline\alpha_n}}
\left\|
\iota_{\frac12}\bigl(C_n(y,y)\bigr)
\right\|_{q,\widehat\phi,\frac12}^{1/2}.
\]
Finally, let $n\to\infty$.  Using
\cref{eq:Lp-conv-convex} and
$\overline\alpha_n\to\alpha$, we obtain
\[
\left\|\iota_{\frac12}\bigl(x-E(x)\bigr)
\right\|_{p,\phi,\frac12}
\le
\frac{p}{\sqrt{2\alpha}}
\left\|
\iota_{\frac12}\bigl(\widehat\Gamma(y,y)\bigr)
\right\|_{q,\widehat\phi,\frac12}^{1/2}.
\]

\medskip
\textbf{Acknowledgments.}
MJ was partially supported by NSF Grant DMS-2247114. JW acknowledges support from the Spring 2025 Graduate Travel Award. The authors thank the Institute for Pure and Applied Mathematics (IPAM), which is supported by the
National Science Foundation (Grant No.~DMS-1925919), for its hospitality, during which the initial
idea of this project was formulated. The authors thank Joel Tropp for helpful discussions.
JW also thanks Roy Araiza, Jihong Cai, David Jekel, Marc Klinger, Rolando de Santiago,
Vincent Villalobos, and Peixue Wu for helpful discussions.

\textbf{Data Availability.}{
Data sharing does not apply to this article, as no datasets were generated or
analyzed during the current study.}

\textbf{Conflict of interest.}{
The authors have no conflicts of interest to declare that are relevant to the
content of this article.}
\bibliography{poincare}
\end{document}